\documentclass[11pt]{article}
\usepackage[utf8]{inputenc}
\usepackage{graphicx,tikz, verbatim, amsmath,subfig,amssymb}

\usepackage{hyperref} 
\usepackage{cleveref}
\usepackage{autonum}

\usepackage[russian,english]{babel}
\usepackage[outdir=./fig/]{epstopdf}
\usepackage[labelfont=bf,justification=centering]{caption}
\usepackage[top=50pt, bottom=50pt, left=50pt, right=50pt]{geometry}
\usepackage{lipsum}
\usepackage{tikz-3dplot,calc}

\usetikzlibrary{patterns,arrows,decorations.markings,positioning,calc,shapes}

\tikzset{middlearrow/.style 2 args={
        decoration={markings,
            mark= at position #1 with {\arrow{#2}} ,
        },        postaction={decorate}    }	}

\newtheorem{lemma}{Lemma}
\newtheorem{remark}[lemma]{Remark}
\newcommand{\dint}{{\rm d}}

\allowdisplaybreaks
\graphicspath{{./fig/}}

    \crefname{figure}{Figure}{Figures}
    \Crefname{figure}{Figure}{Figures}
    \crefname{table}{Table}{Tables}
    \Crefname{table}{Table}{Tables}
    \crefname{section}{\S}{\S}
    \Crefname{section}{\S}{\S}
    \crefname{equation}{}{}
    \Crefname{equation}{}{}
    \crefname{remark}{Remark}{Remarks}
    \Crefname{remark}{Remark}{Remarks}

\title{On the eigenmodes of periodic orbits for \\multiple scattering problems in 2D}

\author{
Daan Huybrechs \\ daan.huybrechs@cs.kuleuven.be\\
Department of Computer Science\\
KU Leuven, Belgium
\and
Peter Opsomer (corresponding author) \\ peter.opsomer@cs.kuleuven.be \\
Department of Computer Science\\
KU Leuven, Belgium
}

\begin{document}

\maketitle


\begin{abstract}

Wave propagation and acoustic scattering problems require vast computational resources to be solved accurately at high frequencies. Asymptotic methods can make this cost potentially frequency independent by explicitly extracting the oscillatory properties of the solution. However, the high-frequency wave pattern becomes very complicated in the presence of multiple scattering obstacles. We consider a boundary integral equation formulation of the Helmholtz equation in two dimensions involving several obstacles, for which ray tracing schemes have been previously proposed. The existing analysis of ray tracing schemes focuses on periodic orbits between a subset of the obstacles. One observes that the densities on each of the obstacles converge to an equilibrium after a few iterations. In this paper we present an asymptotic approximation of the phases of those densities in equilibrium, in the form of a Taylor series. The densities represent a full cycle of reflections in a periodic orbit. We initially exploit symmetry in the case of two circular scatterers, but also provide an explicit algorithm for an arbitrary number of general 2D obstacles. The coefficients, as well as the time to compute them, are independent of the wavenumber and of the incident wave. The results may be used to accelerate ray tracing schemes after a small number of initial iterations.

\end{abstract}


\section{Introduction} \label{Sintro}

Numerical simulations in acoustics are often based on a boundary integral equation reformulation of the Helmholtz equation, see for example \cite{ColtonKress,BEA}. An incoming wave that is scattered by an obstacle $\Omega$ with boundary $\Gamma$ results in a scattered field $u^{\text{s}}(\mathbf{x})$, that can be represented by the so-called \emph{single layer potential},
\begin{equation}
 u^{\text{s}}(\mathbf{x}) = \int_{\Gamma} K(\mathbf{x},\mathbf{y}) v(\mathbf{y}) \dint s(\mathbf{y}). \label{Eslpot}
\end{equation}
Here, $K(\mathbf{x},\mathbf{y})$ is the Green's function of the Helmholtz equation with wavenumber $k$ and $v(\mathbf{y})$ is the unknown \emph{density function}, defined on $\Gamma$. Sound-soft and sound-hard obstacles give rise to Dirichlet or Neumann boundary conditions respectively, with zero pressure or zero normal velocity on $\Gamma$. These lead to integral equations of the first or second kinds. Though other representations of the scattered field exist, which give rise to different integral equations, in this paper we present our results using the representation above. The analytical results depend solely on the parameters of the physical problem, such as the geometry and relative locations of the scattering obstacles, and not on the particular integral equation that is used.

The numerical discretisation of, for example, the Dirichlet problem $u^{\text{s}}(\mathbf{x}) = f(\mathbf{x})$, $\mathbf{x} \in \Gamma$, leads to a large and dense linear system \cite{babich1991shortwave}. A family of hybrid numerical-asymptotic methods aims to significantly reduce the size of the linear system by incorporating information about the solution from asymptotic analysis (see e.g. \cite{abboud,bruno,localSol,unifPhaseExtr,groth}, as well as the review \cite{reviewLangdon} and references therein). In particular, phase-extraction methods use information about the phase $g$ (or multiple phases $g_m$) of the solution $v(\mathbf{y})$, in order to discretize only the remaining non-oscillatory parts $f_m$ in the following factorization:
\[
 v(\mathbf{y}) = \sum_{m=1}^M f_m(\mathbf{y},k)e^{ikg_m(\mathbf{y})}.
\]

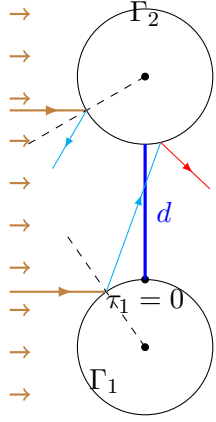
\begin{figure}
\centering
\begin{tikzpicture}[scale = 0.9] 
	\draw (0,0) [fill=white] circle [radius=1];
	\draw (0,4) [fill=white] circle [radius=1];
	\draw (0,0) [fill=black] circle [radius=0.05];
	\draw (-0.6, -0.5) node{$\Gamma_1$};
	\draw (0,4) [fill=black] circle [radius=0.05];
	\draw (0, 4.6) node[above]{$\Gamma_2$};	

	\draw (0,1) [fill=black] circle [radius=0.05] node[below]{$\tau_1=0$};
	
	\def\thone{50}
	\coordinate (c) at ({cos(\thone)}, {sin(\thone)} );
	\def\thtwo{-78.76};
	\coordinate (q) at ({cos(\thtwo)}, {(4+sin(\thtwo))} );
	\coordinate (w) at ({cos(\thtwo)/2}, {(4+sin(\thtwo)/2)} );
	\coordinate (x) at ({cos(\thtwo)*0.05}, {(4+sin(\thtwo)*0.05)} );
	
	\draw[very thick,blue] (0,3) -- (0,1) node[midway, right] (dper) {$d$};
    \foreach \pos in {0,1,...,9}
    {
    	\draw[->,brown, thick] (-2,{-0.7+\pos*5/8}) -- (-1.7,{-0.7+\pos*5/8});
    }
	\draw[brown,thick,middlearrow={0.65}{latex}] (-2,3.5) -- (-0.866, 3.5);
	\draw[dashed] (0,4) -- (-1.732, 3);
	\draw[cyan,middlearrow={0.65}{latex}] (-0.866, 3.5) -- (-1.366, 2.64); 
	
	\draw[brown,thick,middlearrow={0.65}{latex}] (-2, 0.819) -- (-0.574, 0.819); 
	\draw[dashed] (0,0) -- (-1.147,1.638);
	\draw[cyan,middlearrow={0.65}{latex}] (-0.574, 0.819) -- (0.2299, 3.0268); 
	\draw[red,middlearrow={0.65}{latex}] (0.2299, 3.0268) -- (0.9564, 2.3398); 
\end{tikzpicture}
\caption[Reflections in multiple scattering]{A symmetric multiple scattering configuration with two disks. Shown are an incident plane wave (brown), its reflections off the first and second obstacle (cyan), and one further second reflection (red). These rays leave the scene. The only possible rays that can bounce around for many reflections, or even indefinitely, are those close to or exactly on the periodic orbit (blue). Geometrically, the periodic orbit minimizes the distance between the two obstacles.}
\label{FtwoCircRefl}
\end{figure}

The phases $g_m$ can be determined from asymptotic analysis, for example via ray tracing. Ray tracing methods for scattering configurations with multiple obstacles at high frequencies lead to a sequence of single-scattering problems, mimicking exactly the tracing of wavefronts, in such a way that the phase can be extracted from a single density on a single obstacle in each case \cite{geuzaine2005multiple,2DEcevit,3DEcevit,chandlerwilde2015nonconvex}. The methodology is illustrated in \cref{FtwoCircRefl} for the case of two circular obstacles. The first computation is the scattering of the incoming wave by just one of the obstacles, disregarding the other one. In this case, the density mostly inherits the phase of the incoming wave, evaluated along the boundary. Next, the reflected field is considered as an incoming wave for the other obstacle, and so on. This process is repeated and care has to be taken to include all reflections by all obstacles until convergence.

\begin{figure}
\centering 
\tdplotsetmaincoords{60}{120}
\begin{tikzpicture}[scale=1.0,tdplot_main_coords,axis/.style={->,blue,thick},vector/.style={-stealth,red,very thick},vector guide/.style={dashed,red,thick}] 
    \def\a{-2}
    \def\b{3.5}
    \def\rad{0.8}
    \def\c{-5}
    \def\d{-1}
    \def\r{0.6}

    \def\th{60}
    \def\ph{150}
    \def\thpl{-110}
    \coordinate (O) at (\c,\d,0);
    \foreach \angle in {-90,-60,...,90}
    {
        \tdplotsinandcos{\sintheta}{\costheta}{\angle}
        \coordinate (P) at (\c,\d,{\r*\sintheta});
        \tdplotdrawarc{(P)}{{\r*\costheta}}{0}{360}{}{}
        \tdplotsetthetaplanecoords{\angle}
        \tdplotdrawarc[tdplot_rotated_coords]{(O)}{\r}{0}{360}{}{}
    }

    \tdplotsetcoord{P}{1}{\th}{\ph}
    \def\thh{75}
    \def\phh{-90}
    \coordinate (B) at ({\a+(\rad*sin(\thh)*cos(\phh))}, {\b+(\rad*sin(\thh)*sin(\phh))}, {\rad*cos(\thh)} );
    \draw[line width=1mm,blue] (-0.8,0.5,0.1)--({\a+0.2*\rad},{\b-0.7*\rad},{0.3*\rad});
    \draw[line width=1mm,blue] ({\a+0.2*\rad},{\b-0.7*\rad},{0.3*\rad})--({\c+0.7*\r},{\d+0.45*\r},{-0.15*\r});
    \draw[line width=1mm,blue] ({\c+0.7*\r},{\d+0.45*\r},{-0.15*\r})--(-0.8,0.5,0.1);
    \def\fact{1/(sqrt(\a*\a+\b*\b))}
    \def\factR{(sqrt(\a*\a+\b*\b)-\rad)/(sqrt(\a*\a+\b*\b))}
    \draw[line width=1mm,red] ({0+\fact*\a},{\fact*\b},0)--({\a*\factR},{\b*\factR},0);
    \def\fact{1/(sqrt(\c*\c+\d*\d))}
    \def\factR{(sqrt(\c*\c+\d*\d)-\r)/(sqrt(\c*\c+\d*\d))}
    \draw[line width=1mm,red] ({0+\fact*\c},{\fact*\d},0)--({\c*\factR},{\d*\factR},0);
    \def\fact{1/(sqrt((\a-\c)*(\a-\c)+(\b-\d)*(\b-\d) ) ) }
    \draw[line width=1mm,red] ({\a+\rad*\fact*(\c-\a)},{\b+\fact*\rad*(\d-\b)},0)--({\c+\fact*\r*(\a-\c)},{\d+\fact*\r*(\b-\d)},0);

    \coordinate (O) at (0,0,0);
    \foreach \angle in {-90,-60,...,90}
    {
        \tdplotsinandcos{\sintheta}{\costheta}{\angle}
        \coordinate (P) at (0,0,\sintheta);
        \tdplotdrawarc{(P)}{\costheta}{0}{360}{}{}
        \tdplotsetthetaplanecoords{\angle}
        \tdplotdrawarc[tdplot_rotated_coords]{(O)}{1}{0}{360}{}{}
    }
    \coordinate (O) at (\a,\b,0);
    \foreach \angle in {-90,-60,...,90}
    {
        \tdplotsinandcos{\sintheta}{\costheta}{\angle}
        \coordinate (P) at (\a,\b,{\rad*\sintheta});
        \tdplotdrawarc{(P)}{{\rad*\costheta}}{0}{360}{}{}
        \tdplotsetthetaplanecoords{\angle}
	\tdplotdrawarc[tdplot_rotated_coords]{(O)}{\rad}{0}{360}{}{}
    }
\end{tikzpicture}
\caption[Periodic orbits between three spheres]{Depiction of periodic orbits between three spherical scatterers: one periodic orbit visits all three scatterers in a triangle-shaped trajectory minimizing the total path length (blue), and one periodic orbit exists between each combination of two spheres (red).}
\label{FMult}
\end{figure}
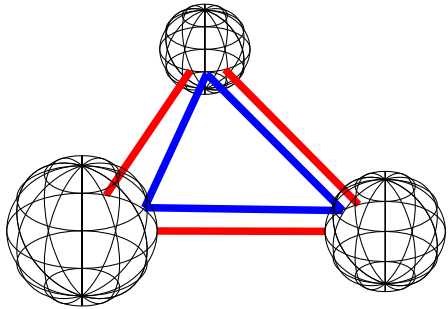

It was shown by Ecevit et al. in a series of papers that for a set of convex obstacles the sequence of single-scattering problems can be organised according to certain \emph{periodic orbits} \cite{2DEcevit,3DEcevit, EcevitOzen2017GBEMConvex,EcevitReitich2005Rate,boubendirEcevit2017acceleration}. Examples of these are illustrated in \cref{FMult}. Each periodic orbit consists of the shortest path between a subset of obstacles, and one can intuitively understand them as follows. Rays originate in a source and reflect off obstacles in the scene. Rays that do not reflect onto another obstacle leave the scene forever. For rays that do reflect again, one can repeat the reasoning on the next iteration of reflections. Fewer and fewer rays will remain, and those that do can only be close to the aforementioned periodic orbits. Indeed, a ray that travels exactly \emph{on} the shortest path between a collection of obstacles remains trapped in the scene forever, while any other ray that deviates from these particular trajectories ultimately leaves the scene.

The rays in periodic orbits settle down to an asymptotic limit -- asymptotic both in the number of iterations and in the frequency -- which we call a \emph{mode} of the configuration of scatterers. From experiments, one observes that the phases of the densities involved converge to some equilibrium after a number of iterations. Indeed, this can be seen in for example \cite[Fig. 1 \& 2]{geuzaine2005multiple}, although this feature was not further investigated there.

Say a periodic orbit involves $J$ obstacles $\Omega_j$, $j=1,\ldots,J$. Rays that travel along the orbit induce a density $V_j$ on each of the participating obstacles. Starting with density $V_1$ on the first obstacle, the observation is that after a full cycle of reflections the induced density on the first obstacle is (asymptotically) proportional to $V_1$ again. We will assume the densities to have the form
\begin{equation}
  V_j(\tau) = v_j^\text{smooth}(\tau) e^{ik\phi_j(\tau)}, \qquad j=1,\ldots,J, \label{EeigvExp} 
\end{equation}
where the phase $\phi_j(\tau)$ corresponds to the conjectured limiting phase described above. We expect that this phase is only well-defined for $\tau$ corresponding to the portions visible to the previous scatterer in the periodic orbit, but $v_j^\text{smooth}(\tau)$ is very small for other $\tau$ anyway. 

The `modes' can be seen as eigenfunctions of the multiple scattering problem. Indeed, the existence of a limit indicates that a wave with the density $V_1$ on the first obstacle induces a density $V_2$ on the second obstacle, $V_3$ on the third obstacle and so on, such that the last density $V_J$ induces precisely the density $V_1$ again on the first obstacle, up to a constant factor. This constant factor can be seen as an eigenvalue, and the collection of densities on each of the obstacles is an eigenfunction of the sequence of single scattering problems. It only depends on the geometry of the scatterers and not on the incident wave.

Similar periodic or \emph{short orbits} have appeared in other settings, in connection to Random Matrix Theory \cite{smilansky,shortOrb}. Certain quantum graphs are employed to study disordered systems. In principle, rays can keep bouncing from an obstacle to another arbitrary obstacle and this wave chaos leads to the application of Random Matrix Theory to obtain spectral statistics of those graphs. Although wave chaos theory has been developed primarily in the context of quantum mechanics, the similarities with statistical methods such as Statistical Energy Analysis are explained in \cite{sondergaard}. The case of an array of cylinders has been studied in \cite{LintonEvans} and can be viewed in analogy with the case of two disks studied here. 

The goal of the current paper is to compute those limiting phases, assuming their existence. We formulate a highly nonlinear system of differential equations that describes the phases. We show that our goal can be achieved if we restrict ourselves to the computation of Taylor series expansions of the phases around the points where the periodic orbits touch the obstacles. This computation is independent of any incident wave, and it is also computationally independent of the wavenumber $k$ of the overall Helmholtz problem and allows solving the system numerically in a simple setting. Still, having an approximation close to the periodic orbit may already be sufficient for practical puposes, as the modulus of $V_j(\tau)$ is maximal for $\tau$ corresponding to rays near the periodic orbit. 
Physically, this is due to rays rapidly leaving the scene the further away they are from the periodic orbit.

Our implementation of the computation of Taylor series and other results in this paper is publicly available on GitHub \cite{github} and a copy of the code can always be obtained from the authors.

Firstly, we exploit symmetry for the simple case of two circular scatterers, as shown in \cref{FtwoCircles}. In this case, which was also presented in \cite{twoCircles}, there is a single periodic orbit, namely the shortest path between the two disks. Using symmetry, it is sufficient to study a single phase function. Secondly, a general procedure is presented to approximate the phase of densities associated with periodic orbits with arbitrary numbers of obstacles including 
quite general geometries. We analyze the limiting phase via its Taylor series around the point where the shortest path hits the boundary. The advantage of the Taylor series approximation is that it leads to fully explicit expressions.

We formulate the problem as an eigenvalue problem for an oscillatory integral operator that represents scattering. Looking for eigenfunctions of the form \cref{EeigvExp} leads to oscillatory integrals. In order to analyze these asymptotically for a large wavenumber $k$, we use the well-known method of steepest descent \cite{Wong,BleisteinHandelsman}. It deforms the path of integration into paths of steepest descent in the complex plane, such that oscillating integrands are transformed to rapidly decaying integrands. A considerable simplification is that we do not need to determine these paths exactly, as we are only interested in the phase of these oscillatory integrals. We do remark that a more extensive analysis of the integrals could lead to an asymptotic description of the amplitudes $v_j^\text{smooth}(\tau)$ in \cref{EeigvExp} as well as the phase.

A similar methodology was employed in \cite{samDaan} in order to obtain the asymptotic expansion of the solution of a scattering problem, possibly involving multiple obstacles as well. In particular, the method of steepest descent was used in \cite{samDaan} to track the forward propagation of rays, given an initial incoming wave. The main difference is that in our current setting, we are concerned with an eigenvalue problem instead.

The paper is structured as follows. First, we detail the formulation of a non-linear system of equations that describes the phases on two obstacles, later generalised to multiple obstacles. In \cref{SserTwoCirc}, we derive explicit Taylor expansions for the phases of the mode of a periodic orbit with two disks. In the next section this is again generalized to a larger number of obstacles with general shapes. A geometrical interpretation of the phase is given in \cref{Sgeometric}. The next section details the numerical approximation of the phase that we use to validate our results in \cref{Sresults}. Finally, we present some concluding remarks.

\begin{remark}\label{rem:conditions}
In this paper we compute the asymptotic phases of the densities in periodic orbits, under the assumption that these orbits exist. This is the case at least when the configuration of obstacles satisfies the conditions stated by Ecevit et al. in \cite{2DEcevit}, where 2D multiple scattering problems involving multiple obstacles were first analyzed. Thus, following \cite{2DEcevit}, we assume a set of $J$ disjoint, smooth and convex obstacles, where in addition rays between two obstacles never intersect a third obstacle (the no-occlusion condition). Yet, we show that the expressions do seem to be valid for some more general configurations in \cref{Sresults}.
\end{remark}

\section{Problem statement}

In this section, we arrive at the formulation of a non-linear system of equations that describes the phases on each obstacle. We focus on $J=2$ scatterers and generalise this to an arbitrary number at the end of this section.

\subsection{Notation}

As in \cref{EeigvExp}, we will refer to functions on obstacle $\Omega_j$ with subscript $j$, where $j=1,\ldots,J$. Each obstacle $\Omega_j$ has a boundary $\Gamma_j$, and we assume a parameterization which we denote by $\mathbf{y} = \Gamma_j(\tau)$. Here, $\mathbf{y}$ is a point on the obstacle corresponding to a unique value $\tau \in [0,1]$ in the parameter domain, which we assume to be the interval $[0,1]$ for all obstacles. We will use $\mathbf{y}$ and $\tau$ interchangeably to denote a point. For clarity, we will often use an indexed parameter variable $\tau_j$ to denote a point on $\Gamma_j$.

We assume that the periodic orbit with $J$ obstacles is known: i.e., the obstacles, their order and the corresponding reflection points $\Gamma_j(\tau_j^*)$ on them. Given an ordered collection of obstacles $\{ \Omega_j \}_{j=1}^J$, the points ${\mathcal T} = \{\tau_j^*\}$ in the corresponding periodic orbit minimise the total length
\begin{equation} \label{eq:length}
	L({\mathcal T}) = \Vert \Gamma_J(\tau_J^*) - \Gamma_{1}(\tau_{1}^*) \Vert + \sum_{j=1}^{J-1} \Vert \Gamma_j(\tau_j^*) - \Gamma_{j+1}(\tau_{j+1}^*) \Vert =  \sum_{j=1}^{J} \Vert \Gamma_j(\tau_j^*) - \Gamma_{j+1}(\tau_{j+1}^*) \Vert.
\end{equation}
In the latter sum as well as further on, we define $\Gamma_{J+1}=\Gamma_1$ in order to shorten the formulae.

We assume that there is a unique solution to this minimisation problem, which is clear if all the obstacles involved are convex \cite{2DEcevit}. However, we do present an example later on with numerical results where this condition is not met. It appears, for example, that the formulae continue to hold if all objects are locally convex and a local minimizer of the form above can be found.

In the expressions above, we have assumed for simplicity of presentation that the periodic orbit contains exactly the obstacles $\Omega_j$, $j=1,\ldots,J$, in that order. Of course, one may consider a periodic orbit between any subset of obstacles, say, $\Omega_3$ and $\Omega_5$. The formulae that follow remain valid simply by considering the multiple scattering problem involving only $\Omega_3$ and $\Omega_5$, i.e., with $J=2$.

The \emph{mode} that corresponds to the orbit ${\mathcal T}$ is the collection of densities $\{ V_j(\tau) \}_{j=1}^J$ on each of the obstacles in ${\mathcal T}$. The distance function between points on two consecutive obstacles is briefly written as
\begin{equation}\label{Edelta}
	\Delta_{i,j}(\tau_i,\tau_j) = \Vert \Gamma_i(\tau_i) - \Gamma_j(\tau_j) \Vert, \qquad 1 \leq i, j \leq J.
\end{equation}

Finally, we call the function $g(x)$ in an asymptotic expression of the form
\[
 f(x; k) \sim u(x; k) e^{ik g(x)}, \qquad k \gg 1,
\]
the \emph{phase} of $f$. The notation above indicates that $f$ and $u$ depend on a wavenumber $k$. We note right away that the phase is not uniquely defined. In the decomposition we assume that $u$ is non-oscillatory, in the sense that its derivatives do not grow with increasing $k$. Even then, any constant $c$ added to $g$ can readily be absorbed into $u$. Hence, $g$ can only be defined up to a constant shift.

\subsection{Coupling operators and the eigenvalue problem}

In order to formulate the eigenvalue problem, recall the single layer potential \cref{Eslpot}. We shall use the notation
\begin{equation}
	(S_{ji} V_i)(\mathbf{x}) = \int_{\Gamma_i} K(\mathbf{x},\mathbf{y}) V_i(\mathbf{y}) \dint s(\mathbf{y}), \qquad \mathbf{x \in \Gamma_j}. \label{EslpSV}
\end{equation}
for its application to a density $V_i$ defined on $\Gamma_i$, and with the point $\mathbf{x}$ on $\Gamma_j$. An incident wave on $\Gamma_j$ is written as $u_j^{\textrm{inc}}$. Disregarding the possibility of resonances, for the time being, the single layer potential leads to coupled integral equations of the first kind for a multiple scattering sound-soft Dirichlet problem. For a set of two obstacles, the coupled system is given by
\begin{equation}
	\begin{cases}  S_{11} m_1 + S_{12} m_2 &= -u_1^{\textrm{inc}}, \\
	S_{21} m_1 + S_{22} m_2 &= -u_2^{\textrm{inc}}. \end{cases} \label{EcouplIE}
\end{equation}
Here, $m_1$ and $m_2$ are the densities on obstacles $\Gamma_1$ and $\Gamma_2$ and the operators $S_{ji}$ are the single layer potentials \cref{EslpSV}. One way to solve this problem, alluded to as ray tracing in the introduction, is to first ignore the other obstacle by leaving out the cross-terms and then to trace consecutive reflections. We numerically illustrate that this converges at the start of \cref{SperOrbitResTwo} and this technique is applied analytically in \cite{samDaan} as well.

The operator under investigation in this paper for $J=2$ is given succinctly by
\begin{equation}
	T =  S_{11}^{-1} S_{12} S_{22}^{-1} S_{21} , \label{ET}
\end{equation}
which represents a full cycle of reflections over two obstacles, independently of the incident wave. With this formulation, the goal of the paper is to find an eigenfunction of $T$ of the form \cref{EeigvExp}, such that
\begin{equation}
 T V_1 = \lambda V_1.
\end{equation}

We proceed by analyzing the operator $T$ asymptotically for large $k$. Expression \cref{ET} consists of a number of steps. It is sufficient to analyze only $S_{11}$ and $S_{21}$, since $S_{22}$ and $S_{12}$ are entirely analogous in the case of just two obstacles. We proceed with the latter operator $S_{21}$ first.

\subsection{Field scattered by the first obstacle onto the second} \label{SperOrbSP}

When we calculate the field that $\Gamma_1$ scatters onto $\Gamma_2$, we obtain
\begin{align}
	S_{21} V_1 & = \int_{\Gamma_1} \frac{i}{4} H_0^{(1)}(k\Delta_{1,2}(\tau, \tau_2)) V_1(\tau) d \tau,
\end{align}
where we have substituted the Green's function $K(x,y)$ of \cref{Eslpot} by its explicit expression for the Dirichlet 2D Helmholtz problem. In order to analyze the integral asymptotically for large $k$, we want to explicitly characterize its total phase.

Since the points $\Gamma_1(\tau_1)$ and $\Gamma_2(\tau)$ are well separated for all $\tau \in [0,1]$, and since $k$ is assumed to be large, we can replace the Hankel function by the leading order term of its expansion for large arguments \cite[(10.17.5)]{DLMF},
\[
 H_0^{(1)}(z) \sim \frac{\sqrt{2}}{\sqrt{\pi z}} e^{i (z - \pi/4)}, \qquad z \to \infty.
\]

Furthermore, we assume $V_1$ to have the form given by \cref{EeigvExp}. Finally, we parameterize the boundary $\Gamma_1$ using the parameterization $\Gamma_1(\tau)$ with the same name and with Jacobian $ \Vert \nabla \Gamma_1(\tau) \Vert$. This results in the asymptotic approximation
\begin{align}
	S_{21} V_1 & \sim \int_0^1 \frac{i}{4} (k\Delta_{1,2}(\tau_1,\tau_2)\pi/2)^{-1/2} e^{ik\Delta_{1,2}(\tau_1,\tau_2)-i\pi/4} v_1^{\text{smooth}}(\tau_1) e^{ik\phi_1(\tau_1)} \Vert \nabla \Gamma_1(\tau_1) \Vert d\tau_1, \label{EexpCircleFar}
\end{align}
where $\Delta_{1,2}$ is the Euclidean distance between two points on the two obstacles, recall \cref{Edelta},
\[
\Delta_{1,2}(\tau_1,\tau_2) = \Vert \Gamma_2(\tau_2) - \Gamma_1(\tau_1) \Vert.
\]

We are mainly interested in the leading order behaviour of the integral. Higher order terms can be found by substituting the full expansion of the Hankel function, a process detailed in \cite{samDaan}. Yet, this entails also a full asymptotic expansion of $v_1^{\text{smooth}}(\tau_1)$, which as of yet is unknown. Importantly, the total phase of the integrand \cref{EexpCircleFar} is fully explicit,
\begin{equation}\label{eq:g21}
 g_{2,1}(\tau_2,\tau_1) = \Delta_{1,2}(\tau_1,\tau_2) + \phi_1(\tau_1).
\end{equation}

As is standard in asymptotic expansions of integrals \cite{Wong, anCont, BleisteinHandelsman}, one notes that the leading order contribution to an oscillatory integral originates 
in a region around the stationary points of the phase of the integrand. A stationary point $\chi$ of a phase function $g(x)$ is such that $g'(\chi)=0$, i.e., it is a root of the derivative of $g$. In the context of scattering, such stationary points correspond precisely to points of reflection \cite{samDaan}.

In order to further approximate \cref{EexpCircleFar} asymptotically, we look for a stationary point $\chi_1(\tau_2)$ for each value of $\tau_2$. Physically, this means that a ray hitting obstacle $\Gamma_2$ at the point $\tau_2$ originated in the point $\tau_1 = \chi_1(\tau_2)$ on $\Gamma_1$: it is like inverse ray tracing. The condition for $\chi_1(\tau_2)$ is that the derivative of the total phase vanishes with respect to the integration variable $\tau$, hence
\begin{equation}
	\frac{\partial g_{2,1}}{\partial \tau_1}(\tau_2, \chi_1(\tau_2)) = 0. \label{EcircleSP}
\end{equation}
This is an implicit definition of the function $\chi_1$. The phase of the resulting function $S_{12} V_1$ is given by
\begin{equation}
 g_{2,1}(\tau_2, \chi_1(\tau_2)) = \Delta_{1,2}(\chi_1(\tau_2), \tau_2) + \phi_1(\chi_1(\tau_2)).
\end{equation}
It equals the phase $\phi_1$ at the reflection point, plus the distance from that point to $\Gamma_2(\tau_2)$. This agrees with the principle that the phase along a ray equals the distance travelled.

Recall that $\tau_j^*$ are the points of the periodic orbit. We know by construction that $\chi_1(\tau_2^*) = \tau_1^*$: this expresses that the periodic orbit contains the ray from $\tau_1^*$ to $\tau_2^*$. Thus, though $\chi_1$ is characterized only implicitly above, we know at least its value in one point. Furthermore, due to the conditions on the obstacles, the function is necessarily at least locally single-valued.

Finally, we can also define the function $\chi_2(\tau_1)$ to find the point on $\Gamma_2$ from which a ray reflects onto $\Gamma_1(\tau_1)$. Note that $\chi_2$ is not the inverse of $\chi_1$, but we do have in this scenario of just two obstacles in the orbit that both
\begin{equation}\label{eq:initial}
\chi_1(\tau_2^*) = \tau_1^* \quad \mbox{and} \quad \chi_2(\tau_1^*) = \tau_2^*.
\end{equation}

\subsection{Single scattering problem on the first obstacle} \label{SperOrbitGS}

Next, we consider $S_{11}^{-1}$. A full computation of the asymptotic expansion of the solution to a single scattering problem, even with phase extraction, is rather involved. The computation is carried out in \cite{samDaan}. The problem of computing the phase of the solution is simpler since, in fact, it is known that the phase of the solution equals the phase of the incoming wave from $\Gamma_2$ under the conditions in \cite{2DEcevit}. As such, 
\begin{equation} 
 \phi_1(\tau_1) = g_{1,2}(\tau_1, \chi_2(\tau_1)) = \Delta_{2,1}(\chi_2(\tau_1), \tau_1) + \phi_2(\chi_2(\tau_1)). \label{EcirclePhase}
\end{equation}

For completeness, the process of obtaining full asymptotics is as follows. Assuming a density of the form \cref{EeigvExp}, the application of the integral operator on $\Gamma_1$ yields
\begin{equation}\label{eq:S11_V1}
 S_{11} V_1(\tau_1) = u(\tau_1) = \int_{\Gamma_1} \frac{i}{4} H_0^{(1)}\left(k \Vert \Gamma_1(\tau_1) - \Gamma_1(\tau)\Vert \right) \left(v_1^{\text{smooth}}(\tau) e^{ik\phi_1(\tau)} \right) \Vert \nabla \Gamma_1(\tau) \Vert d\tau. 
\end{equation}
A technical complication is the appearance of the logarithmically singular Hankel function, which in this case can not be replaced by its expansion for large arguments as before in \cref{EexpCircleFar}, where it led to a simpler oscillatory exponential. An alternative is to substitute the Hankel function by its Mehler-Sonine integral representation, which does have a complex exponential in the integrand \cite[(10.9.10)]{DLMF}
\begin{equation}
       H_0^{(1)} (z) = \frac{1}{\pi i} \int_{-\infty-\pi i/2}^{+\infty+\pi i/2} e^{i z \cosh(t)} \dint t. \label{EintHankel}
\end{equation}
After splitting the integration interval around $\tau_1^*$, inverting the distance function and interchanging integration variables, a regular steepest descent analysis can be applied on the resulting bivariate integral. 

The main observation to make here is that the only relevant asymptotic contribution to the singular oscillatory integral originates in the singularity at $\tau_1=\tau$. As one can apply a path of steepest descent such that the Mehler-Sonine integral representation of $H_0^{(1)}$ does not contribute to the phase of the integrand (because the singularity is precisely where the argument of $H_0^{(1)}$ is zero), the phase of $S_{11} V_1$ indeed equals the phase of $V_1$.

\subsection{System of nonlinear equations for two obstacles}

We have sufficient information to describe the action of the operator $T$ given by \cref{ET} on functions of the form \cref{EeigvExp}. We look for a function $V_1$ with phase $\phi_1(\tau)$, such that the density after a full cycle of reflections has the phase $\phi_1(\tau) + \mu$, where $\mu$ represents a constant phase shift. Using \cref{EcircleSP}--\cref{EcirclePhase}, the nonlinear system of differential-algebraic equations that $\phi_1(\tau)$, $\phi_2(\tau)$, $\chi_1(\tau)$ and $\chi_2(\tau)$ should satisfy is:
\begin{equation}
\begin{cases}
	\frac{\partial g_{2,1}}{\partial \tau_1}(\tau_2, \chi_1(\tau_2)) = 0 \\ 
	\phi_2(\tau_2) = \Delta_{1,2}(\chi_1(\tau_2),\tau_2)+\phi_1(\chi_1(\tau_2))  \\
	\frac{ \partial g_{1,2}}{\partial \tau_2}(\tau_1, \chi_2(\tau_1)) = 0 \\
	\phi_1(\tau_1) = \Delta_{2,1}(\chi_2(\tau_1),\tau_1)+\phi_2(\chi_2(\tau_1)) - \mu.  
\end{cases} \label{Esys2Obst}
\end{equation}
Analogously to \eqref{eq:g21}, we have introduced the phase of the integral representation of $S_{12}V_2$,
\[
 g_{1,2}(\tau_1, \tau_2) = \Delta_{2,1}(\tau_2,\tau_1) + \phi_2(\tau_2).
\]

The first three equations in \eqref{Esys2Obst} were introduced before. The final equation introduces the constant phase shift $\mu$, and ensures that $\phi_1$ is recovered after a full orbit up to $\mu$. Since the phase equals the distance travelled, this implies that
\[
 \mu = L({\mathcal T}),
\]
where $L({\mathcal T})$ is the length of the periodic orbit \cref{eq:length}.

We have thus derived at a system of $4$ equations, which is highly nonlinear in the four unknown functions $\phi_1(\tau)$, $\phi_2(\tau)$, $\chi_1(\tau)$ and $\chi_2(\tau)$.  Even in the case of two circular obstacles with equal radius, the system does not seem to allow an explicit analytic solution.

It is not obvious whether the system has a solution at all and, if so, whether the solution is unique. 

\begin{remark}
 We make no rigorous claims regarding the (unique) solvability of \cref{Esys2Obst} globally for $\tau_1,\tau_2 \in [0,1]$. Indeed, since our extraction of phases is based on geometrical optics, it is not valid in the shadow region of the obstacles (the shadow relative to the previous obstacle in the orbit). By the no-occlusion condition of \cite{2DEcevit} and formulated in Remark \ref{rem:conditions}, the periodic orbit is such that rays close to it are amenable to geometrical optics. We assert that a solution to \cref{Esys2Obst} exists in a neighbourhood of the points $\tau_j^*$ on each $\Gamma_j$, $j=1,\ldots,J$.
\end{remark}

In order to find a local solution to \cref{Esys2Obst}, we augment the system with the initial conditions \cref{eq:initial}. Furthermore, it is clear that the phase $\phi_1$ is only determined up to a constant, which we may choose, e.g., by prescribing a value for $\phi_1(\tau_1^*)$.

In the following section, we will resort to Taylor series expansions of the unknown functions to approximate the local solution of \cref{Esys2Obst}. This approach does lead to fully explicit expressions for the coefficients. The series is computed around the periodic orbit, where the mode $V_j(\tau)$ has a (nearby) local maximum. Though our derivation is not presented in a fully mathematically rigorous way in this paper, the fact that the computation succeeds is a strong indication that the system \cref{Esys2Obst}, along with the initial conditions for $\chi$ and the function value $\phi_1(\tau^*)$, is indeed uniquely solvable in a neighourhood of the periodic orbit.

\subsection{System for a general number of obstacles}

We can readily generalize the arguments of this section to an arbitrary number of obstacles $J$. The orbit operator $T$ becomes
\begin{equation} 
 T = S_{11}^{-1} S_{1J} \ldots S_{33}^{-1} S_{32} S_{22}^{-1} S_{21}.\label{ETJ}
\end{equation}

We can consider each pair of obstacles involved in the reflection of rays, the obstacles $(j+1,j)$, for $j=1,\ldots,J-1$, and the final pair $(J,1)$. An increasing $j$ means the rays progress forward in time through the obstacles and hence diverge away from the periodic orbit. The function $\chi_j(\tau_{j+1})$ instead contracts around the periodic orbit, i.e., it traces the rays backwards in time from obstacle $j+1$ back to $j$.

We define the phase function of the integral representation of $S_{ji} V_i$ (recall \eqref{eq:S11_V1}) as
\begin{equation}\label{eq:gji}
  g_{j,i}(\tau_j, \tau_i) = \Delta_{i,j}(\tau_i, \tau_j) + \phi_i(\tau_i),
\end{equation}
where the distance $\Delta_{i,j}$ between points on obstacles $\Gamma_i$ and $\Gamma_j$ is given by \cref{Edelta}.

With that notation, we can formulate the following $2J \times 2J$ nonlinear system:
\begin{align} \label{EsysGenObst}
\begin{cases}
  \frac{\partial g_{j+1,j}}{\partial \tau_j }(\tau_{j+1}, \chi_j(\tau_{j+1})) = 0,    \\
	 \phi_{j+1}(\tau_{j+1}) = \Delta_{j,j+1}(\chi_j(\tau_{j+1}),\tau_{j+1}) + \phi_j(\chi_j(\tau_{j+1})),\\  
	 \frac{\partial g_{1,J}}{\partial \tau_{J} }(\tau_1, \chi_J(\tau_1)) = 0,    \\ 
	 \phi_1(\tau_1) = \Delta_{J,1}(\chi_{J}(\tau_1),\tau_1)+\phi_J(\chi_{J}(\tau_1)) - \mu.
\end{cases} \qquad j = 1,\ldots,J-1. 
\end{align}
Here, as before, the phase shift after a full cycle equals the length of the periodic orbit,
\[
 \mu = L({\mathcal T}).
\]
This means that the phase of the ray that follows the periodic orbit increases exactly by the length of the path that is travelled.

In order to characterize a local solution to \cref{EsysGenObst}, we augment the system with the initial conditions
\[
 \chi_j(\tau_{j+1}^*) = \tau_j^*,
\]
and we choose a value for $\phi_1(\tau_1^*)$.

\section{The special case of two disks} \label{SserTwoCirc}

In this section, we derive explicit Taylor expansions for the phases $\phi_1$ and $\phi_2$ of the mode $\{ V_1, V_2 \}$ of a periodic orbit between two disks.

\subsection{Parameterization of the boundary}

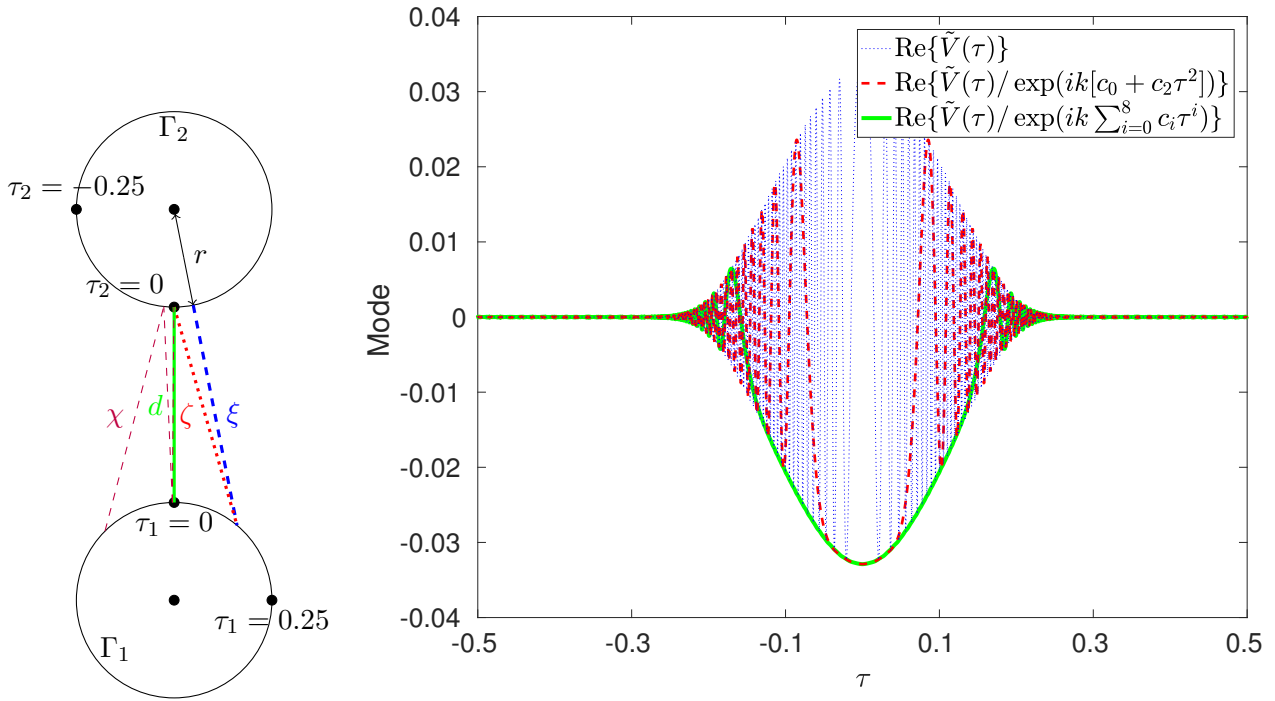
\begin{figure}[ht]
\centering
\subfloat{
\begin{tikzpicture}[scale = 1.3]
	\draw (0,0) [fill=white] circle [radius=1];
	\draw (0,4) [fill=white] circle [radius=1];
	\draw (0,0) [fill=black] circle [radius=0.05];
	\draw (-0.6, -0.5) node{$\Gamma_1$};
	\draw (0,4) [fill=black] circle [radius=0.05];
	\draw (0, 4.6) node[above]{$\Gamma_2$};	

	\draw (0,3) [fill=black] circle [radius=0.05] node[above left]{$\tau_2=0$};
	\draw (-1,4) [fill=black] circle [radius=0.05] node[above]{$\tau_2=-0.25$};
	\draw (0,1) [fill=black] circle [radius=0.05] node[below]{$\tau_1=0$};
	\draw (1,0) [fill=black] circle [radius=0.05] node[below]{$\tau_1=0.25$};
	
	\def\thone{50}
	\coordinate (c) at ({cos(\thone)}, {sin(\thone)} );
	\def\thtwo{-78.76};
	\coordinate (q) at ({cos(\thtwo)}, {(4+sin(\thtwo))} );
	\coordinate (w) at ({cos(\thtwo)/2}, {(4+sin(\thtwo)/2)} );
	\coordinate (x) at ({cos(\thtwo)*0.05}, {(4+sin(\thtwo)*0.05)} );
	
	\draw[<-] (q) -- (w) node[right] {$r$};
	\draw[<-] (x) -- (w);
	\draw[very thick,blue,dashed] (q) -- (c)  node[midway, right] {$\xi$};
	\draw[very thick,red,dotted] (0,3) -- (c) node[midway, left] {$\zeta$};

	\draw[very thick,green] (0,3) -- (0,1) node[midway, left] {$d$};

	\def\tone{-45}
	\coordinate (c1) at ({sin(\tone)}, {cos(\tone)} );
	\def\ttwo{-6.0769}
	\coordinate (c2) at ({sin(\ttwo)}, {4-cos(\ttwo)} );
	\def\tthree{-1.0376}
	\coordinate (c3) at ({sin(\tthree)}, {cos(\tthree)} );
	\def\tfour{-0.1780}
	\coordinate (c4) at ({sin(\tfour)}, {4-cos(\tfour)} );
	\draw[purple,dashed] (c1) -- (c2) node[midway,left] {$\chi$};
	\draw[purple,dashed] (c3) -- (c2);
	\draw[purple,dashed] (c3) -- (c4);
\end{tikzpicture} 
}
\subfloat{\includegraphics[width=0.73\hsize]{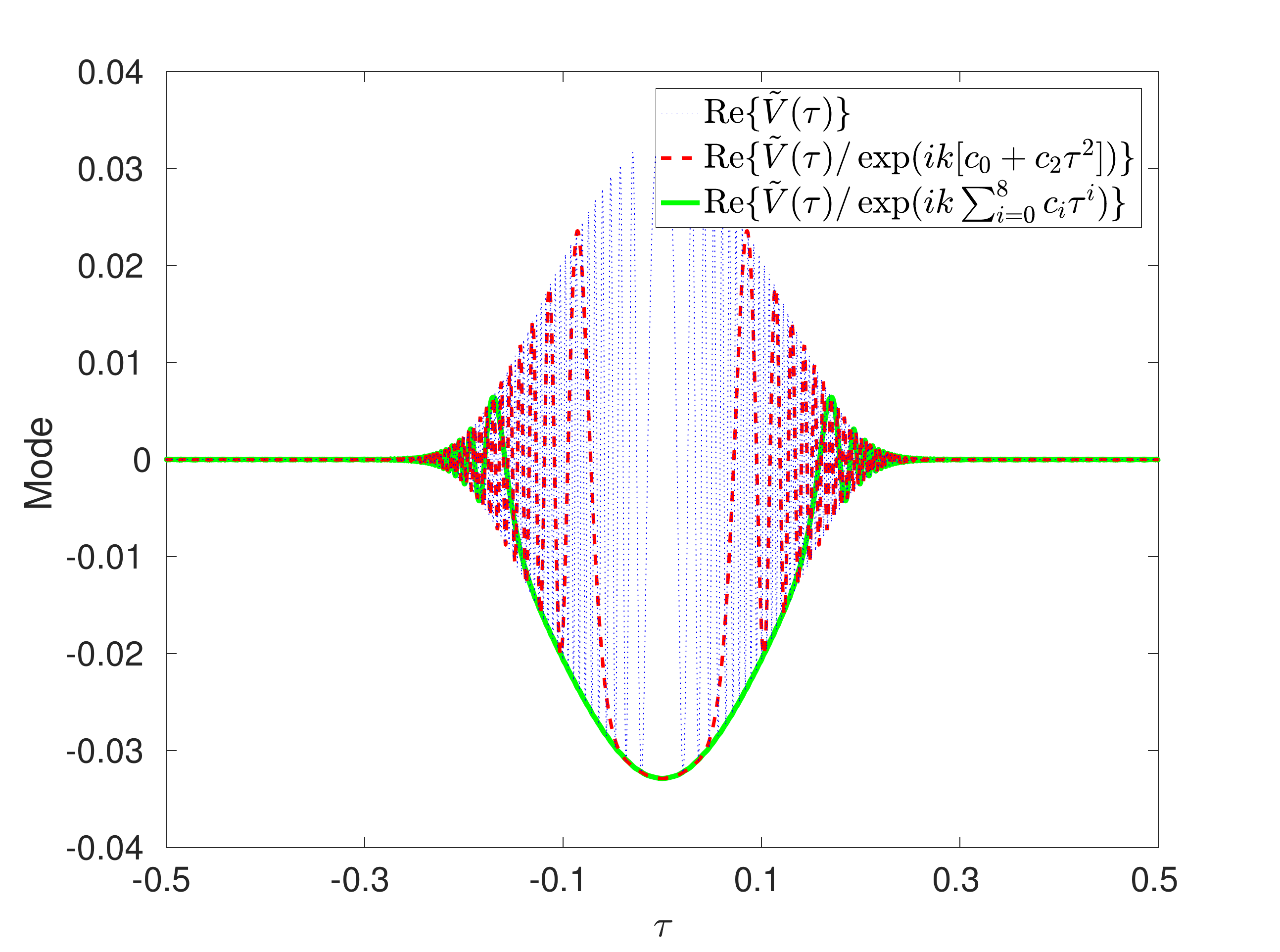} }
\caption[Two disks with eigenvector]{Two disks under consideration in this paper (left) and a numerical illustration of the mode $\tilde{V}(\tau)$ via the first eigenvector $E_{j,1}$ of the discrete operator $M$ (see \cref{Snumerical}) with different phases factored out (right). The wavenumber $k=2^9$. The right panel shows that extracting a phase using its Taylor series expansion of degree $2$ or degree $8$ from the numerical solution leads to a remainder that is non-oscillatory near $\tau^*=0$.}
\label{FtwoCircles}
\end{figure}

The non-overlapping disks in \cref{FtwoCircles} are separated by a distance $d$ and the corresponding shortest path is the periodic orbit along which rays will be trapped forever. We choose to parameterize the circular boundaries with radius $r$ specifically as, with $\tau \in [0,1]$,
\begin{align}
	\Gamma_{1}(\tau) & = r [\sin(2\pi\tau), \cos(2\pi\tau)], \label{EparCircOne} \\
	\Gamma_{2}(\tau) & = [0,d+2r] + r[\sin(2\pi\tau), -\cos(2\pi\tau)].
\end{align} 
Since the length of the periodic orbit is twice the distance, we have $\mu = L({\mathcal T}) = 2d$. The points in the periodic orbit are $\tau_1^* = \tau_2^*=0$. When we show numerical values further on, we choose the parameters 
\[
 r = 1/2 \qquad \mbox{and} \qquad d=1.
\]

The unknown phase function for the density on $\Gamma_1$ was denoted by $\phi_1(\tau)$ before. Due to the geometrical symmetry, and because the parameterization is defined in opposite directions on the two obstacles, we have that $\phi_2(\tau) = \phi_1(\tau) + \nu$, where $\nu$ is a constant phase factor. It follows from \cref{Esys2Obst} at $\tau^*=0$ that $\nu = d$, hence
\[
 \phi_2(\tau) = \phi_1(\tau)+d.
\]
In the following we omit the subscript, i.e. $\phi_1(\tau) = \phi(\tau)$. Similarly, we write $\chi_1(\tau) = \chi_2(\tau) = \chi(\tau)$.

\subsection{Procedure} \label{SprocTwoCirc} 

We look for a Taylor series around $\tau^* = 0$ of the form 
\begin{equation}
  \phi(\tau) \sim \sum_{i=0}^\infty c_i (\tau - \tau^*)^i = \sum_{i=0}^\infty c_i \tau^i. 
\end{equation}
We are free to choose the constant $c_0$, since it corresponds purely to a constant phase shift and eigenfunctions are only determined up to a constant factor. We arbitrarily set $c_0=d$.

We also expand $\chi(\tau)$ into a Taylor series,
\begin{align}
	\chi(\tau) & \sim \sum_{i=0}^\infty a_i \tau^i, 
\end{align}
as well as the residuals of the first two equations in \eqref{Esys2Obst}:
\begin{align}
	\frac{\partial g_{21}}{\partial \tau_1}(\tau, \chi(\tau)) & \sim \sum_{i=0}^\infty \psi_i \tau^i. \label{EpsiTwoCirc} \\ 
	\phi(\chi(\tau)) + d - \Delta_{1,2}(\tau,\chi(\tau)) - \phi(\tau) & \sim \sum_{i=0}^\infty \omega_i \tau^i. \label{EomegaTwoCirc} 
\end{align}
Note that the third and fourth equation of \eqref{Esys2Obst} lead to the exact same expansion, recalling that $\phi_2(\tau) = \phi_1(\tau) + d$ and $\mu = 2d$.

We already know that $a_0 = 0$ since $\chi(\tau_1^*) = \tau_2^* = 0$. Substituting $\tau=0$ into \eqref{EomegaTwoCirc} also shows that $\omega_0 = 0$. All other coefficients are computed in a recursive procedure in the next subsection. For each index $i$, starting from $i=0$, we solve
\begin{equation}\label{eq:step}
	\omega_{i+1} = 0 \qquad \mbox{and} \qquad \psi_i = 0, \qquad i=0,1,\ldots,
\end{equation}
for the unknowns $a_i$ and $c_{i+1}$. These conditions express that the residuals of the nonlinear equations vanish to all orders. The newly found coefficients for index $i$ can be used in the next iteration for index $i+1$.

\subsection{Explicit Taylor series coefficients}\label{Sexplicit}

In the first iteration of \eqref{eq:step} for $i=0$, starting from $\omega_{1} = 0 = \psi_0$ we find that
\begin{equation}
	a_1 c_1 - c_1 = 0 = c_1 \quad \Rightarrow \quad c_1=0.
\end{equation}
This corresponds to $\phi'(0) = 0$, which indicates that $\phi$ has a local extremum at the point $\tau^* = 0$ as we expect.

For $i=1$, we need the series expansion of the left hand side of \cref{EomegaTwoCirc} up to order two, and of the left hand side of \cref{EpsiTwoCirc} up to order $1$. Substituting the expansions for $\chi$ and $\phi$ leads to
\begin{align} 
	\mathcal{O}(\tau^3) & = c_0 + c_1 (a_0 + a_1\tau + a_2 \tau^2) + c_2 (a_0 + a_1\tau + a_2 \tau^2)^2  - c_0 - c_1\tau -c^2\tau^2 \\ 
	& + \sqrt{ \left(\frac{\sin 2\pi\tau}{2} - \frac{\sin \left(2\pi (a_0 + a_1\tau + a_2 \tau^2)  \right) }{2}\right)^2 + \left(\frac{\cos 2\pi\tau}{2} -2 + \frac{\cos \left( 2\pi (a_0 + a_1\tau + a_2 \tau^2) \right) }{2} \right)^2}, \\ 
	\mathcal{O}(\tau^2) & = \left. \frac{\partial \left\{   \sqrt{ \left( \frac{\sin 2\pi\tau}{2} - \frac{\sin 2\pi \tau_\alpha }{2} \right)^2 + \left( \frac{\cos 2\pi\tau}{2} -2 + \frac{\cos 2\pi\tau_\alpha}{2} \right)^2 }    + c_0 + c_1\tau_\alpha + c^2\tau_\alpha^2 \right\} }{\partial \tau_\alpha} \right. ,
\end{align}
where the latter is evaluated at $\tau_\alpha = a_0 + a_1\tau + a_2 \tau^2$. We also expand the distance function $\Delta_{1,2}$ into a Taylor series. This requires the expansion of the sines, cosines and the square root function in the expressions above, which can all be performed explicitly. Inserting the newly found value $c_1=0$, the equations \eqref{eq:step} become the system of quadratic equations
\begin{align}
	\omega_2 & = 0 =\frac{3}{2}\pi^2 a_1^2 - \pi^2 a_1 + a_1^2 c_2 + \frac{3\pi^2}{2} - c_2, \label{EtwoCircOm2} \\
	\psi_1 & = 0 = 3\pi^2 a_1 - \pi^2 + 2a_1 c_2, \label{EtwoCircPsi2}
\end{align}
from which we deduce
\begin{equation}
	c_2 = \sqrt{2}\pi^2, \quad a_1=3-2\sqrt{2}.
\end{equation}

Because \eqref{EtwoCircOm2} is quadratic in $a_1$, there is a second solution, which is given by $a_1=3+2\sqrt{2}, c_2 = -\sqrt{2}\pi^2$. This solution is not physically relevant, as the left side of \cref{Fphase} indicates that $c_2 > 0$ and, furthermore, $a_1 > 1$ would imply that the rays move away from the periodic orbit shown in the left part of \cref{FtwoCircles}. Instead, we want $\chi$ to be contracting.

For higher $i$, the systems of equations that arise are linear in all coefficients and the solution is unambiguous. Due to the symmetry in the current problem of two disks, we find that $a_i=c_{i+1}=0$ for even values of $i$. Further computation leads to the next nonzero coefficients
\begin{equation}
	c_4 = \frac{-11}{12}\sqrt{2}\pi^4, \qquad a_3 = -7\pi^2(17\sqrt{2}-24), 
\end{equation}
and
\begin{align}
	c_6 & =\frac{2783\sqrt{2}\pi^6}{2520}, & & a_5 = - \frac{\pi^4}{84}(1205811\sqrt{2} - 1705312), \\ 
	c_8 & = \frac{-358021}{205632} \sqrt{2}\pi^8, & & a_7 = \frac{-\pi^6}{128520}(289615597399\sqrt{2} - 409578202752).
\end{align}
Recall that these values are specific to the choice $d=1$ and $r=1/2$, though they can easily be generalized by repeating the derivation.

The accuracy of these formulae is depicted graphically in the right panel of \cref{FtwoCircles}. We have computed a numerical approximation to the eigenfunction of the operator $T$ (described further on in \cref{Snumerical}), and extracted the Taylor series approximation of the phase. The mode is non-oscillatory in a neighbourhood of the critical point $\tau=0$, it is supported in the region where the disks can `see' each other, and the oscillatoriness decreases with increasing number of terms in the expansion.

\section{The case of multiple obstacles with general convex shapes}

We derive explicit Taylor expansions for the phases of the mode $\{ V_j \}_{j=1}^J$ of a periodic orbit with a general number $J$ of obstacles with a general shape. The process is largely comparable to the special case of two disks, yet slightly more involved. One issue that complicates the algebra is the expansion of the distance function $\Delta_{j,j+1}(\tau_j,\tau_{j+1})$, which is a function of two arguments. It could be done explicitly for the case of two disks, but in general it can only be expanded in terms of the Taylor series expansions of the parameterization of the boundaries at the points $\tau_j^*$.

\subsection{Computation of the periodic orbit}

In the analytical derivation that follows, we assume that the periodic orbit ${\mathcal T} = \{ \tau_j^* \}_{j=1}^J$ is known, given the ordered collection of $J$ obstacles. In practice, of course, one has to compute the reflection points $\Gamma_j(\tau_j^*)$ on each of the obstacles $\Gamma_j$. Recall that the periodic orbit minimizes the sum of the distances between each pair of consecutive obstacles. To that end, one has to solve the minimization problem \cref{eq:length}.

In our implementation, we find the corresponding parameters $\tau_j^*$ on each obstacle $\Gamma_j$ by a nonlinear optimisation, more specifically the quasi-Newton algorithm in \texttt{fminunc} in \textsc{Matlab} based on supplied gradient and Hessian, and with all tolerances set to machine precision. The initialisation is via the parameters $\tau_j$ giving the point approximately closest to the mean of $100$ approximately equidistant points on each obstacle. Though the solution of the minimisation problem is unique for a collection of convex obstacles, the numerical procedure typically also converges to a periodic orbit for non-convex obstacles. We present an example further on in \cref{SresThreeObsts}.

\subsection{Series expansion of the bivariate distance function} \label{SperOrbitSettGenObst} 

The distance function $\Delta_{j,j+1}(\tau_j,\tau_{j+1})$ appears in the system \cref{EsysGenObst} and has to be expanded in both of its arguments. For this purpose, we determine a double series expansion of the distance between two points $\tau_{j}$ and $\tau_{j+1}$ on the obstacles $\Gamma_j$ and $\Gamma_{j+1}$, where we define $\Gamma_{J+1}=\Gamma_1$ for convenience of notation.

Since we work in two dimensions, the parameterization $\Gamma_j$ is vector-valued. We define its vector-valued Taylor series in terms of the Taylor series of the $x$ and $y$-components,
\begin{equation}\label{Egammapar}
	\Gamma_j(\tau_j) = \sum_{n=1}^\infty \begin{pmatrix} \Gamma_{j,x,n} \\ \Gamma_{j,y,n} \end{pmatrix} (\tau_j - \tau_j^*)^{n-1}. 
\end{equation}
The coefficients $\Gamma_{j,x,n}$ and $\Gamma_{j,y,n}$ can be obtained directly from the parameterization $\Gamma_j(\tau_j)$ by differentiation in the point $\tau_j^*$.\footnote{Note that the index $n$ starts at $1$ here, chosen such that there is an exact match between the expressions in this paper and the formulae in our implementation in \textsc{Matlab} \cite{github} (in which vector indices start at $1$).}

Substituting \cref{Egammapar} into expression \cref{Edelta} for the distance function leads, after some algebraic manipulation further detailed below, to
\begin{align}
	\Delta_{j,j+1}(\tau_j, \tau_{j+1}) & = \Vert\Gamma_j(\tau_j) - \Gamma_{j+1}(\tau_{j+1})\Vert \\ 
	& \sim \sqrt{\sum_{l=1}^\infty \sum_{i=1}^\infty (\Lambda_{j,x,i,l} + \Lambda_{j,y,i,l}) (\tau_j - \tau_j^*)^{i-1} (\tau_{j+1} - \tau_{j+1}^*)^{l-1} } \\ 
	& \sim \sqrt{\Lambda_{j,x,1,1} + \Lambda_{j,y,1,1}} \left(1+\sum_{m=1}^\infty {1/2 \choose m} \sum_{k=1}^\infty \sum_{i=\max(1,m-k+2)}^\infty z_{j,m,i,k} (\tau_j -\tau_j^*)^{i-1} (\tau_{j+1} - \tau_{j+1}^*)^{k-1} \right). 
\end{align}
Here, we have introduced several intermediate variables. The values $\Lambda_{j,x,i,l}$ and $\Lambda_{j,y,i,l}$ arise from the full expansion of the square of the difference of two expansions of the form \cref{Egammapar} for $j$ and $j+1$. For the $x$-component, the relation is
\[
 \left( \sum_{m=1}^\infty \Gamma_{j,x,m} (\tau_j - \tau_j^*)^{m-1} - \sum_{l=1}^\infty \Gamma_{j+1,x,l} (\tau_{j+1} - \tau_{j+1}^*)^{l-1}\right)^2 = \sum_{m=1}^\infty \sum_{l=1}^\infty \Lambda_{j,x,m,l} (\tau_j - \tau_j^*)^{m-1} (\tau_{j+1} - \tau_{j+1}^*)^{l-1},
\]
from which we derive the expression for $\Lambda_{j,x,m,l}$ after the calculation of two convolutions:
\begin{equation}\label{Elambda}
	\Lambda_{j,x,m,l} = \delta_{m,1} \left(\sum_{i=1}^l \Gamma_{j, x, i} \Gamma_{j, x, l-i+1} \right)   -2 \Gamma_{j, x, l} \Gamma_{j+1, x, m} + \delta_{l,1} \left(\sum_{i=1}^m \Gamma_{j+1, x, i} \Gamma_{j+1, x, m-i+1} \right). 
\end{equation}
The expression for $\Lambda_{j,y,i,l}$ is entirely analogous. We have used the Kronecker delta ($\delta_{n,k} =1$ if $n=k$ and zero elsewhere) to describe the constant terms when $m=1$ or $l=1$.

The coefficients $z_{j,m,i,k}$ arise from the expansion of the square root, and they are given by
\[
z_{j,1,i,k} = \frac{\Lambda_{j,x,i,k} + \Lambda_{j,y,i,k}}{\Lambda_{j,x,1,1} + \Lambda_{j,y,1,1}},  \qquad \quad  z_{j,m,i,k} = \sum_{r=1}^k \sum_{s=1 + \delta_{r,1} }^{i -\max(0,m-k+r-1)} z_{j,1,s,r} z_{j,m-1, i-s+1, k-r+1}.
\]
The coefficient $z_{j,1,1,1}$ is not used. As an implementation note, when the two last indices in $z_{j,m,i,k}$ are smaller than $a$ the range of the subscript $m$ relevant for the computation of the expansion is limited to $1 \leq m \leq 2a-2$.

The final expression for the distance function has the form
\begin{equation}
	\Delta_{j,j+1}(\tau_j, \tau_{j+1}) \sim \sum_{l=1}^\infty  \sum_{n=1}^\infty f_{j,l,n} (\tau_j - \tau_j^*)^{l-1} (\tau_{j+1} - \tau_{j+1}^*)^{n-1},
\end{equation}
where $f_{j,1,1}$ is the distance between two consecutive points in the orbit,
\[
	f_{j,1,1} = \Vert \Gamma_j(\tau_j^*) -\Gamma_{j+1}(\tau_{j+1}^*) \Vert \quad = \sqrt{\Lambda_{j,x,1,1} + \Lambda_{j,y,1,1} } \quad = \sqrt{\left( \Gamma_{j,x,1} -\Gamma_{j+1,x,1} \right)^2 +\left( \Gamma_{j,y,1} -\Gamma_{j+1,y,1} \right)^2 }.
\]
This is the value around which the square root has been expanded above.

Finally, the coefficients that we are looking for in this section are
\begin{equation}
 f_{j,i,k} = \frac{\partial^i \partial^k \Delta_{j,j+1}}{i! (k!) \partial \tau_j^i \partial \tau_{j+1}^j}(\tau_{j}^*, \tau_{j+1}^*) = \sqrt{\Lambda_{j,x,1,1} + \Lambda_{j,y,1,1} } \sum_{m=1}^{k-2+i} {1/2 \choose m} z_{j,m,i,k}.
\end{equation}
Since $\tau_j^*$ is part of the periodic orbit which minimizes the total distance, we have that
\[
 \frac{\partial \Delta_{j-1,j}}{\partial \tau_j}(\tau_{j-1}^*, \tau_j^*) + \frac{\partial \Delta_{j,j+1}}{\partial \tau_j}(\tau_j^*, \tau_{j+1}^*) = 0.
\]
This leads to the property
\begin{equation}\label{Eminimumproperty}
 f_{j,2,1} = - f_{j-1,1,2}.
\end{equation}

\subsection{Procedure} 

More general than in \cref{SprocTwoCirc}, the points $\tau_j^*$ may be nonzero. We look for a Taylor series of the form
\begin{equation}\label{EgeneralTaylor}
	\phi_j (\tau_j) \sim \sum_{i=0}^\infty c_{j,i} (\tau_j - \tau_j^*)^i. 
\end{equation}
As before, we are free to choose the constant $c_0$ and we set $c_{j,0}=d_j$, the distance to the next obstacle on the periodic orbit. For notational convenience we again define $\Gamma_{J+1} = \Gamma_1$, $\phi_{J+1} = \phi_{1}$ and so on.

We expand $\chi$ and the residuals of the equations in \cref{EsysGenObst} into Taylor series as well: 
\begin{align}
	\chi_j(\tau_{j+1}) & \sim \sum_{i=0}^\infty a_{j,1,i} (\tau_{j+1} - \tau_{j+1}^*)^i, \\
	\phi_{j}(\chi_{j}(\tau_{j+1})) + \Delta_{j,j+1}(\chi_j[\tau_{j+1}],\tau_{j+1}) - \phi_{j+1}(\tau_{j+1}) & \sim \sum_{i=0}^\infty \omega_{j,i} (\tau_{j+1} - \tau_{j+1}^*)^i, \label{Eomega} \\
	\frac{\partial g_{j+1,j}}{\partial \tau_j}(\tau_{j+1}, \chi_j(\tau_{j+1})) & \sim \sum_{i=0}^\infty \psi_{j,i} (\tau_{j+1} - \tau_{j+1}^*)^i. \label{Epsi}
\end{align}

We already know that $a_{j,1,0} = \tau_j^*$ since $\chi_{j}(\tau_{j+1}^*) = \tau_{j}^*$. The coefficients $\omega_{j,i}$ and $\psi_{j,i}$ can be computed explicitly in terms of $a_{j,1,i}$ and $c_{j,i}$, using the bivariate Taylor series expansion of the distance function derived above in \cref{SperOrbitSettGenObst}. All other coefficients are computed in a recursive procedure explained in the next subsection. As before, we want the coefficients in the expansion of the residuals to be zero, with the exception of $\omega_{j,0}$ which may be nonzero because of \cref{EcirclePhase}.

\subsection{Formulae for general 2D obstacles}

We explicitly expand the power of an expansion using a convolution, leading to
\[
  \left[\sum_{n=1}^\infty a_{j,1,n} (\tau_{j+1} - \tau_{j+1}^*)^n \right]^i = \sum_{l=i}^\infty a_{j,i,l} (\tau_{j+1} - \tau_{j+1}^*)^{l},
\]
with
\[
  a_{j,i,l} = \sum_{k=i-1}^{l-1} a_{j,i-1,k} a_{j,1,l-k}
\]
and with $a_{j,i,i} = (a_{j,1,1})^i$.

Then, \cref{Eomega} and \cref{Epsi} respectively become
\begin{align}
	\omega_{j,i} = 0 & = c_{j,i} a_{j,i,i} + \left( \sum_{m=1}^{i-1} c_{j,m} a_{j,m,i} \right)  - c_{j+1,i} + \left(\sum_{r=2}^{i+1} f_{j,r,1} a_{j,r-1,i} \right)  \label{EciOmega}  \\ 
	& + f_{j,1,i+1} + \sum_{n=2}^i \sum_{m=n-1}^{i-1} f_{j,n,i+1-m} a_{j,n-1,m}, \\ 
	\psi_{j,i-1} = 0 & = f_{j,2,i} + \left(\sum_{l=1}^{i-1} \sum_{n=3}^{i-l+2} f_{j,n,l} (n-1) a_{j,n-2,i-l} \right)   +  c_{j,1} \delta_{i,1} \label{EciPsi} \\ 
	& + \left( \sum_{k=2}^{i-1} c_{j,k} k a_{j,k-1,i-1} \right)  + \begin{cases} i c_{j,i} a_{j,i-1,i-1} & i \geq 2, \\ 0 & i = 1. \end{cases} 
\end{align}

This gives a set of $2J$ equations in $2J$ unknowns $a_{j,1,i}$ and $c_{j,i}$. For $i=1$, this becomes
\begin{align}
	c_{j,1} = -f_{j,2,1}, \qquad \qquad c_{j+1,1} = +f_{j,1,2}. \label{Ea1c1} 
\end{align}
Note that, when considering the analogous equations for the obstacles $j+1$ and $j+2$, the first equation above becomes $c_{j+1,1} = -f_{j+1,2,1}$. Combined with the second equation above, this is consistent with \cref{Eminimumproperty}.

\begin{figure}
\centering
\includegraphics[width=0.65\textwidth]{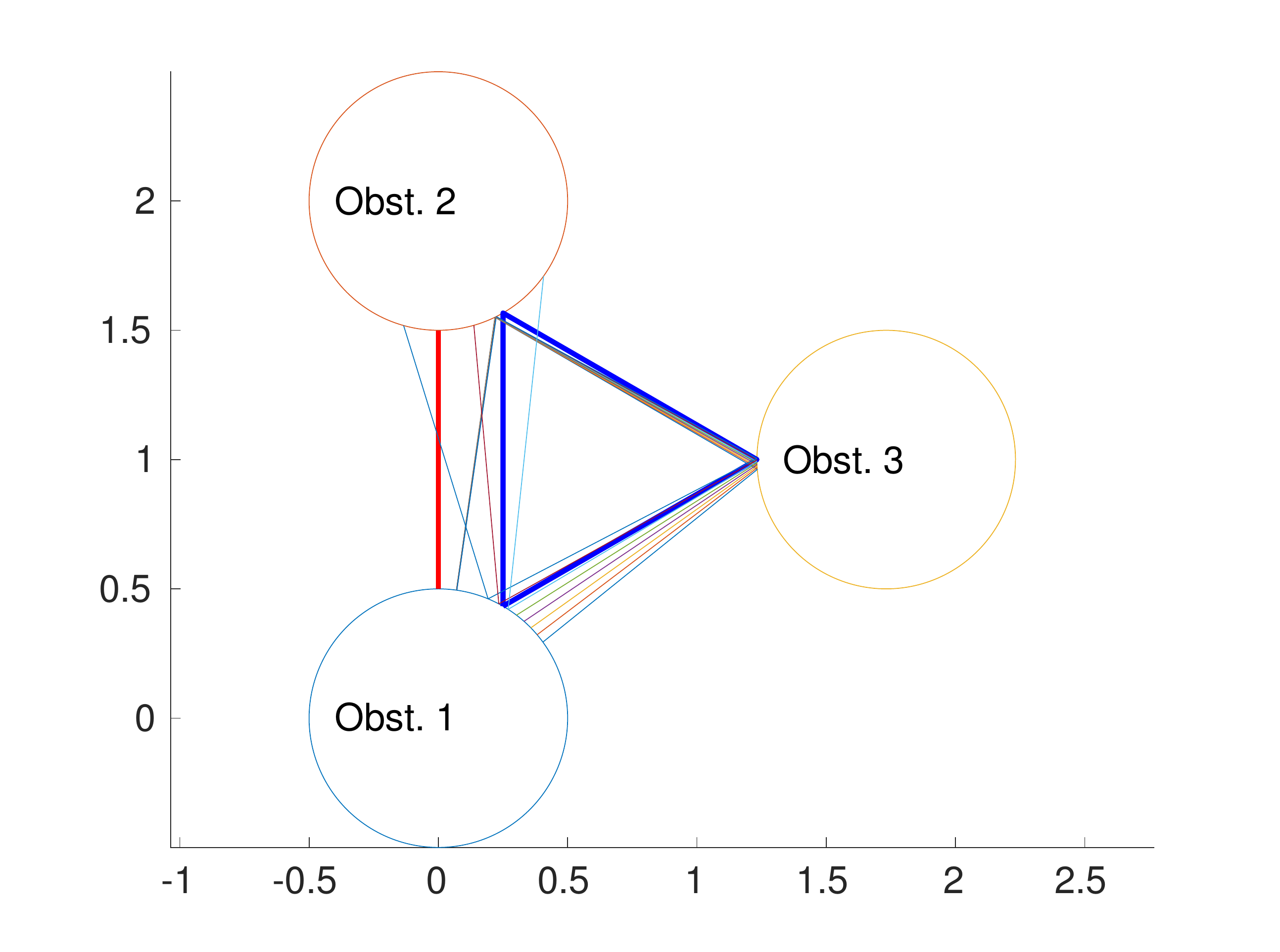}
\caption[Minimum of the phase]{In the case of a periodic orbit between two obstacles (red), the limiting phase on one obstacle has a minimum precisely at the critical point $\tau_j^*$ on the periodic orbit. However, this is not the case for orbits between three or more obstacles (blue). There, the minimum of the phase on obstacle $\Gamma_1$ is the point where a ray leaves orthogonal to $\Gamma_1$ and reaches the periodic orbit after an infinite number of reflections: we have shown 3 or 4 reflections of a bundle of rays starting near the minimum of the phase.}
\label{FminPhPerOrbits}
\end{figure}

From the Taylor series \cref{EgeneralTaylor} it is clear that the value $c_{j,1}$ is the derivative of the phase at the point $\tau_j^*$ on the periodic orbit. If it is zero, the density has an extremum at the critical point, which is related to that the periodic orbit is exactly the shortest distance between the obstacles. Since the signs of $c_{j,1}$ and of $f_{j,2,1}$ are opposite, the phase $\phi_{j}(\tau_j)$ decreases when moving towards the minimum of the distance $\Delta_{j,j+1}(\tau_{j}, \tau_{j+1}^*)$. This minimum does not correspond to the minimum of the phase $\phi_j'(\tau_j^M)=0$. Physically, the latter should be the point where a ray leaves orthogonal to $\Gamma_j$ and reaches the periodic orbit after an infinite number of reflections, where $\chi_j(\tau_{j+1}^*) = \tau_j^* \neq \tau_J^M$ if $c_{j,1} \neq 0$, see \cref{FminPhPerOrbits}. Note that the phase depends on the direction of the periodic orbit as well now, in contrast to a periodic orbit with only two obstacles. 

The coefficient $a_{j,1,1}$ is still unknown. So in the next steps over increasing $i$, all $\omega_{j,i}$ are combined with $\psi_{j,i-1}$ to obtain $c_{j,i}$ and $a_{j,1,i-1}$ simultaneously for all $j$. For $i=2$, \cref{EciPsi} and \cref{EciOmega} become
\begin{align}
	\psi_{j,1} = 0 & = f_{j,2,2} + 2 f_{j,3,1} a_{j,1,1} + 2 c_{j,2} a_{j,1,1}, \\
	\omega_{j,2} = 0 & = [c_{j,2} + f_{j,3,1} ] (a_{j,1,1})^2 + f_{j,2,2} a_{j,1,1}  - c_{j+1,2}  + f_{j,1,3} = f_{j,2,2} a_{j,1,1}/2  - c_{j+1,2}  + f_{j,1,3}.
\end{align}
The latter simplification was possible due to the equation for $\psi_{j,1}$. This gives a system of $2J$ quadratic equations which can be solved with a nonlinear solver with an initial guess $c_{j+1,2} = f_{j,1,3}$ 
and $a_{j,1,1}=0$. A motivation for the latter is that $\chi_{j}(\tau_{j+1})$ should be a function that contracts around the periodic orbit, so $a_{j,1,1}$ should be small: more specifically, $|a_{j,1,1}| < |f_{j,1,2} / f_{j,2,1} |$. We add this check to our code as well as $c_{j,2} > 0$, because we expect $\tau_j^*$ to lie close to a (local) minimum of $\phi_j$ in analogy to the case $J=2$. The initial guess $c_{j+1,2} = f_{j,1,3}$ makes sure we already satisfy $\omega_{j,2}$, one of both equations. 
In the case of the two disks, $f_{1,2,2} = -\pi^2, f_{1,3,1} = 3\pi^2/2 = f_{1,1,3}$ such that the above equations correspond to \cref{EtwoCircOm2} and \cref{EtwoCircPsi2}. We expect that if $f_{j,i,m} = 0$ for all $i+m \leq n+1$ and $j \in [1,J]$, a polynomial system of $2J$ equations, each of which has a total degree $ \leq n$, will give the solution $a_{j,1,1}$ and $c_{j,n}$ for $1\leq j\leq J$. Physically, an increasing $n$ means that both obstacles become flatter at the periodic orbit, resulting in less oscillatory modes.

For $i \geq 3$ in \cref{EciOmega} and \cref{EciPsi}, the following $2J \times 2J$ linear system was implemented
\begin{align}
	\omega_{j,i} = 0 & = c_{j,i}(a_{j,1,1})^i  - c_{j+1,i} + \left( \sum_{j=3}^{i-1} c_{j,m} a_{j,m,i} \right) +  \{c_{j,2} + f_{j,3,1}\} \left( \sum_{k=2}^{i-2} a_{j,1,k} a_{j,1,i-k} \right) \\
	&  \left(\sum_{r=4}^{i+1} f_{j,r,1} a_{j,r-1,i} \right) + f_{j,1,i+1} + \sum_{n=2}^i \sum_{m=n-1}^{i-1 - \delta_{n,2} } f_{j,n,i+1-m} a_{j,n-1,m}, \\ 
	\psi_{j,i-1} = 0 & = i c_{j,i} \{a_{j,1,1}\}^{i-1}  + \left\{2 f_{j,3,1} + 2 c_{j,2} \right\} a_{j,1,i-1} \\ 
	& +  f_{j,2,i} + \left(\sum_{l=1}^{i-1} \sum_{n=3+\delta_{l,1}}^{i-l+2} f_{j,n,l} (n-1) a_{j,n-2,i-l} \right)   +  \sum_{k=3}^{i-1} c_{j,k} k a_{j,k-1,i-1}. 
\end{align}

The validity of our approach is confirmed by the numerical experiments further on involving three obstacles.

\section{A geometrical interpretation of $\chi$ and $\phi$} \label{Sgeometric}

We can assign a geometrical meaning to the stationary point function $\chi(\tau)$ and the phase $\phi(\tau)$ in the case of two disks, by expanding the analogy to ray tracing. Moreover, the reasoning extends to general periodic orbits.

Recall that $\chi(\tau_1)$ is the point on $\Gamma_2$ from which a ray hits the point $\Gamma_1(\tau_1)$. By symmetry, the converse also holds: $\chi(\chi(\tau_1))$ yields a point back on $\Gamma_1$. In the sequence $\chi(\chi(\tau_1)) \to \chi(\tau_1) \to \tau_1$, it should be true that the incoming and reflected ray at the point $\chi(\tau_1)$ on $\Gamma_2$ are at equal angles with the normal direction. In our setting of two disks or radius $r$ and at a distance $d$ of each other, this leads to the equation
\begin{align}
	0 & = 4\pi \chi(\tau_1) + \arctan\left(\frac{r\sin(2\pi\chi(\tau_1)) -r\sin(2\pi\chi(\chi(\tau_1)))}{d+2r -r\cos(2\pi\chi(\tau_1)) -r\cos(2\pi\chi(\chi(\tau_1))) } \right) \label{Echi} \\
	& + \arctan\left( \frac{r\sin(2\pi\chi(\tau_1)) -r\sin(2\pi\tau_1)}{d+2r -r\cos(2\pi\chi(\tau_1)) -r\cos(2\pi \tau_1) } \right). \nonumber
\end{align}
This equation is again highly nonlinear, and it is hard to obtain a symbolic solution on account of the application of $\chi$ on itself. However, one can compute the Taylor series of this expression for $\tau_1$ near zero, and equate all its coefficients to zero, and we checked this to result in the same $a_i$ as in the procedure described before.

One can also numerically approximate $\chi$ using \cref{Echi}. In turn, we can obtain the phase in a limited range of $\tau$ from \eqref{eq:g21}--\eqref{EcircleSP}. Indeed, differentiating $\eqref{eq:g21}$ with respect to $\tau_1$, and substituting into \eqref{EcircleSP}, leads to
\[
 \phi_1'(\chi_1(\tau_2)) = - \frac{\partial \Delta_{1,2}}{\partial \tau_1}(\chi_1(\tau_2),\tau_2).
\]
Multiplying both sides by $\chi_1'(\tau_2)$ and integrating with respect to $\tau_2$ yields the explicit expression
\begin{equation}
 \phi_1(\chi_1(\tau)) = -\int_{\tau_2^*}^\tau \frac{\partial \Delta_{1,2}}{\partial \tau_1}(\chi_1(\tau_2),\tau_2) \chi_1'(\tau_2) \dint \tau_2. \label{EphiChi}
\end{equation}
This determines $\phi_1$, up to inversion of $\chi(\tau)$ in the argument of $\phi_1$ above. The accuracy of the approximation given by \cref{EphiChi} is included in the right panel of \cref{Fphase}, which shows uniform accuracy in a range of $\tau$. This could be extended to a collection of $J$ general obstacles, but that would require repeated integration and differentiation.

We may characterize the phase further in geometrical terms. As in ray tracing, one can track the phase of a ray in a homogeneous medium along straight lines. The phase of a ray is equal to its initial phase plus the travelled distance. Hence, the phase is the sum of the lengths of all the reflecting rays. We define the $r$-times iterated application of $\chi$ by $\chi^{[r]}(\tau)$, such that $\chi^{[0]}(\tau) = \tau$, $\chi^{[1]}(\tau) = \chi(\tau)$, $\chi^{[2]}(\tau) = \chi(\chi(\tau))$ and so on. We say that a point $\Gamma_1(\tau)$ is in the range of the critical point $\tau^*$ if the repeated application of $\chi$ leads to the critical point,
\[
 \lim_{r \to \infty} \chi^{[r]}(\tau) = \tau^*.
\]
This amounts to inverse ray tracing.

Summing the lengths of the reflected rays leads to, for $\tau$ in the range of $\tau^*$,
\begin{equation}
 \phi(\tau) = \phi(\tau^*) + \sum_{r=0}^\infty \left( \Delta_{1,2}\left(\chi^{[2r]}(\tau), \chi^{[2r+1]}(\tau)\right) + \Delta_{2,1}\left(\chi^{[2r+1]}(\tau), \chi^{[2r+2]}(\tau)\right) - 2d\right). \label{EphiInf}
\end{equation}
This formula holds for the case of two disks. Each term in the sum corresponds to a cycle from $\Gamma_1$ to $\Gamma_2$ and back. Note that we subtract the length of the periodic orbit $L({\mathcal T})=2d$ in each cycle, otherwise the sum would diverge. The formula agrees with the repeated application of the fourth and second equation in \cref{Esys2Obst}, each time recursively replacing $\phi_1$ or $\phi_2$ in the right hand sides by their expression in terms of the other function.

For a collection of $J$ obstacles and a periodic orbit ${\mathcal T}$, the reasoning can be extended by tracing a point $\tau_1$ on $\Gamma_1$ back to the critical point $\tau_1^*$ on the periodic orbit. One sums the lengths of each cycle of reflections, following $\chi_J$ down to $\chi_1$, and subtracts in each cycle the total length $L({\mathcal T})$ of the orbit. As above, this corresponds to the recursive application of the second and fourth equation of \cref{EsysGenObst}.

The phase at $\tau_1$ thus computed is the distance of the shortest path traveled by the ray to the periodic orbit, in agreement with the principle of Fermat. At each obstacle, there is reflection at equal angles, and in each cycle the distance between the obstacles along the periodic orbit is subtracted.


\section{Numerical approximation} \label{Snumerical}

We briefly sketch a numerical method to validate the expressions obtained in this paper, based on an implementation of the standard boundary element method for the coupled system \cref{EcouplIE}. Each operator $S_{ji}$ results in a matrix $A_{j,i}$, and we can realize the operator $T$ given by \cref{ET} (for two obstacles) or \cref{ETJ} (for multiple obstacles) numerically.

We assume for each boundary $\Gamma_j$ a set of functions $\psi_{j,n}(\tau)$ to discretize the densities with $N_j$ unknowns, see \cref{EVnumerical} further on. The discretization of a coupled problem with $J$ obstacles results in a dense BEM matrix with block structure form
\begin{equation}
	A = \begin{pmatrix} A_{1,1} & A_{1,2} & \ldots \\ A_{2,1} & A_{2,2} & \ldots \\ \vdots & \vdots & \ddots \end{pmatrix}. \label{EselfCoupling} 
\end{equation}
Here, $A_{j,i} \in \mathbb{C}^{N_j \times N_i}$ represents the discretization of $S_{ji}$. Thus, the blocks $A_{j,j}$ on the diagonal are the BEM matrices of the single scattering problems involving only $\Gamma_j$. The off-diagonal blocks $A_{j,i}$ represent the coupling from $\Gamma_i$ to $\Gamma_j$. Our experiments are based on a collocation approach, in which we have used piecewise linear basis functions and $N_j = 10 k_j$ evenly spaced collocation points on each boundary $\Gamma_j$ \cite{huybrechs2008gratings}.

In the ray-tracing scheme for multiple scattering problems, one only uses the inverse of the diagonal blocks, as these correspond to a single scattering problem~\cite{2DEcevit}. It is shown in \cite{2DEcevit} that, under certain conditions, the norms of the off-diagonal blocks are sufficiently small to allow the iterative process to converge. In addition to the convexity and non-occlusion condition already mentioned, the obstacles should be separated by more than the inverse of the wavenumber. More specifically, regularity estimates are available if $\Delta_{j,i}(\tau_j, \tau_i) \geq 1/k$ for all $1 \leq i,j \leq J$ and $\tau_i,\tau_j \in [0,1]$, see \cite[\S 4.3]{andrew} and \cite{ChGiLaMo:15}.

The numerical analogue of $T$, a full cycle of reflections in a periodic orbit given by \cref{ETJ}, is obtained simply by replacing any operator $S_{ji}$ with its discretization $A_{j,i}$,
\begin{equation}\label{EM}
 M = (A_{1,1}^{-1} A_{1,J}) (A_{J,J}^{-1} A_{J,J-1}) \ldots (A_{2,2}^{-1} A_{2,1}).
\end{equation}
We can numerically compute the eigenvalue decomposition $M= E D E^{-1}$, ordered such that its largest (in modulus) eigenvalue is $D_{1,1}$. Experiments confirm that the complex argument of $D_{1,1}$ equals $L({\mathcal T})$, the length of the periodic orbit, up to a multiple of $2\pi/k$. We proceed with only the first eigenvalue and eigenvector in what follows, noting only that the further eigenvalues seemingly decay rapidly, while the phases of the corresponding eigenvectors are similar to that of the first eigenvector. The modulus of the eigenvalues seems to only depend on the geometry and not on (high) $k$ nor the number of points per wavelength, where increasing the latter two only adds negligible eigenvalues.

We obtain an approximation to the phase of the eigenmode from the phase of \cref{EVnumerical}, where the coefficients in the expansion are given by the first eigenvector. Recall the factorization of \cref{EeigvExp} of $V_j$ into an oscillatory and non-oscillatory part, which is approximated numerically using the first eigenvector in $E$,
\begin{equation} \label{EVnumerical}
	V_j(\tau) = v_j^\text{smooth}(\tau) e^{ik\phi_j(\tau)} \approx \tilde{V}_j(\tau) = \sum_{n=1}^{N_j} E_{n,1} \psi_{j,n}(\tau).
\end{equation}
To obtain the first eigenvectors on other obstacles than the first, one can reorder \cref{EM} and recompute the eigenvalue decomposition $M= E D E^{-1}$, but in our implementation we simply multiply $E_{n,1}$ on the first obstacle with the corresponding pairs of submatrices in \cref{EM}. We can not merely consider the complex phase of point evaluations of $V_j$ to approximate $\phi_j$, since $v_j^\text{smooth}$ may also be complex-valued. However, assuming that the complex argument of $v_j^\text{smooth}$ varies much more slowly than Arg$(\exp[ik\phi_j(\tau)])$ for moderately large $k$, we can estimate the phase $\phi_j$ at a point $\tau$ from a finite difference with small positive or negative $\delta$,
\begin{equation}
	\tilde{\phi}_j(\tau) = \tilde{\phi}_j(\tau-\delta) + \frac{\text{Arg}\left(\tilde{V}_{j}(\tau) \right) - \text{Arg}\left(\tilde{V}_j(\tau-\delta)\right) }{k} - \frac{2\pi}{k} \left\lfloor \frac{\text{Arg}\left(\tilde{V}_{j}(\tau)\right) - \text{Arg}\left(\tilde{V}_j(\tau-\delta)\right) }{2\pi}\right\rfloor.  \label{Ephitildeprox} 
\end{equation}
Note that the final term simply removes any jump of size $2\pi$ in the argument that may have occurred from one point to the next. In our implementation, $|\delta|$ is the distance between collocation points in the $\tau$-domain and because we use piecewise linear basis functions on each boundary, $\tilde{V}_j(\tau-\delta)$ simply equals $E_{n,1}$ at the collocation point $\tau_{1,n} = \tau-\delta$. The numerical phase function was already illustrated in the right panel of \cref{FtwoCircles} and used to validate our Taylor series expressions for $\phi_1(\tau) = \phi_2(\tau)$ in the case of two disks. Several other examples are given in the next section.

Though in this paper we present no analysis of the amplitude function $v_j^\text{smooth}$ of the eigenmode, it seems concentrated in a region around the critical point $\tau_j^*$. As one would expect, the support is effectively limited to that part of the obstacles that `see' each other. The amplitude decays by a factor equal to the modulus of the largest eigenvalue after each cycle of reflections. This is also illustrated in the next section.

\section{Numerical results} \label{Sresults}

\subsection{Two disks} \label{SperOrbitResTwo}

We have already illustrated the accuracy of the Taylor series approximation to the limiting phase function $\phi$ in the right panel of \cref{FtwoCircles}. There, the eigenmode was computed numerically using the numerical method outlined in the previous section, i.e., as the eigenvector of $M$ given by \cref{EM}. The Taylor series of the phase is computed using the explicit expressions for the coefficients given in \cref{Sexplicit}. The figure confirms, at least visually, that the remainder after factoring out the Taylor series of the phase is non-oscillatory, in a neighbourhood of the critical point $\tau^*$.

We elaborate on the case of two disks with further experiments that are more quantitative.

\begin{figure}[h]
\centering
\subfloat{\includegraphics[width=0.45\hsize]{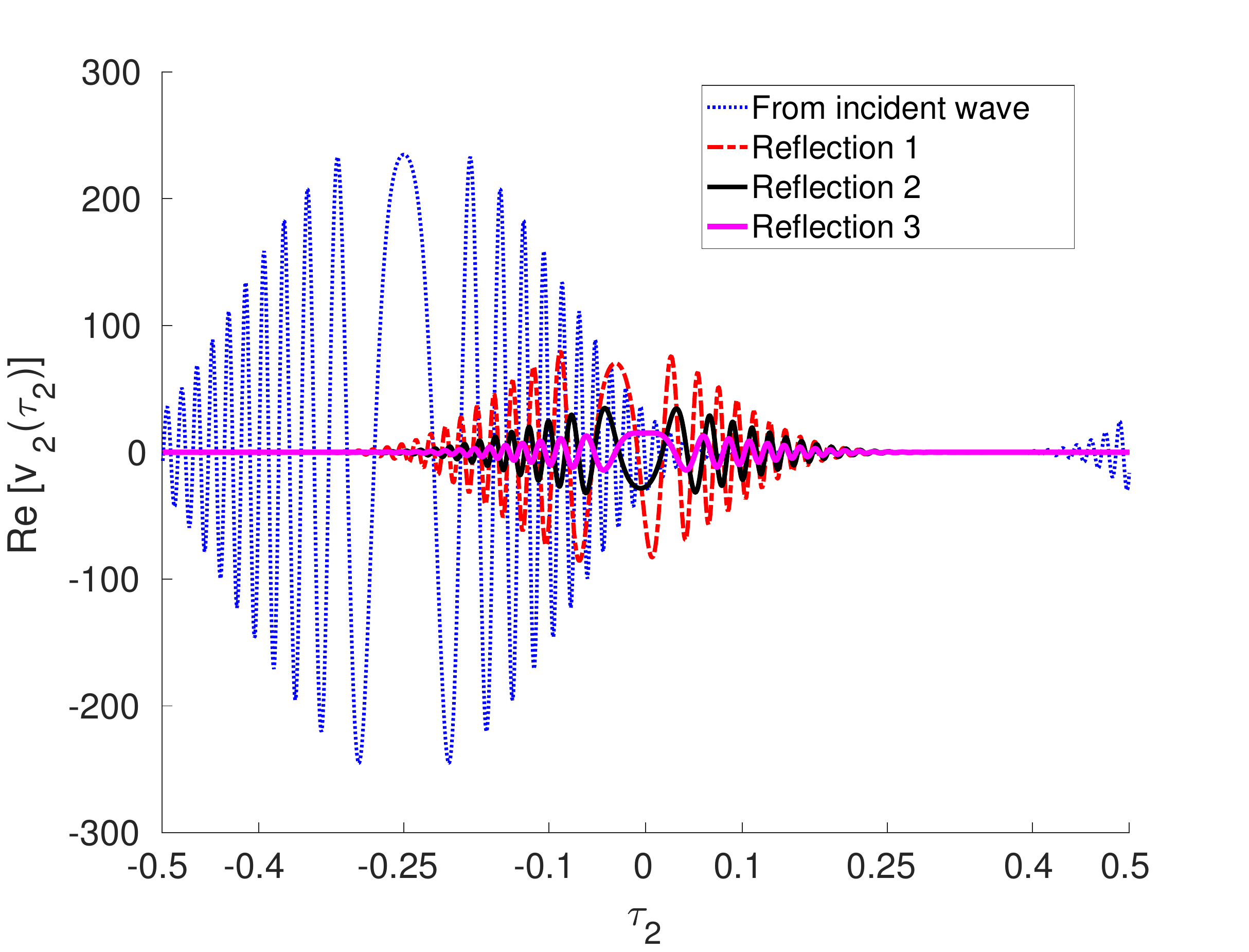} }
\subfloat{\includegraphics[width=0.54\hsize]{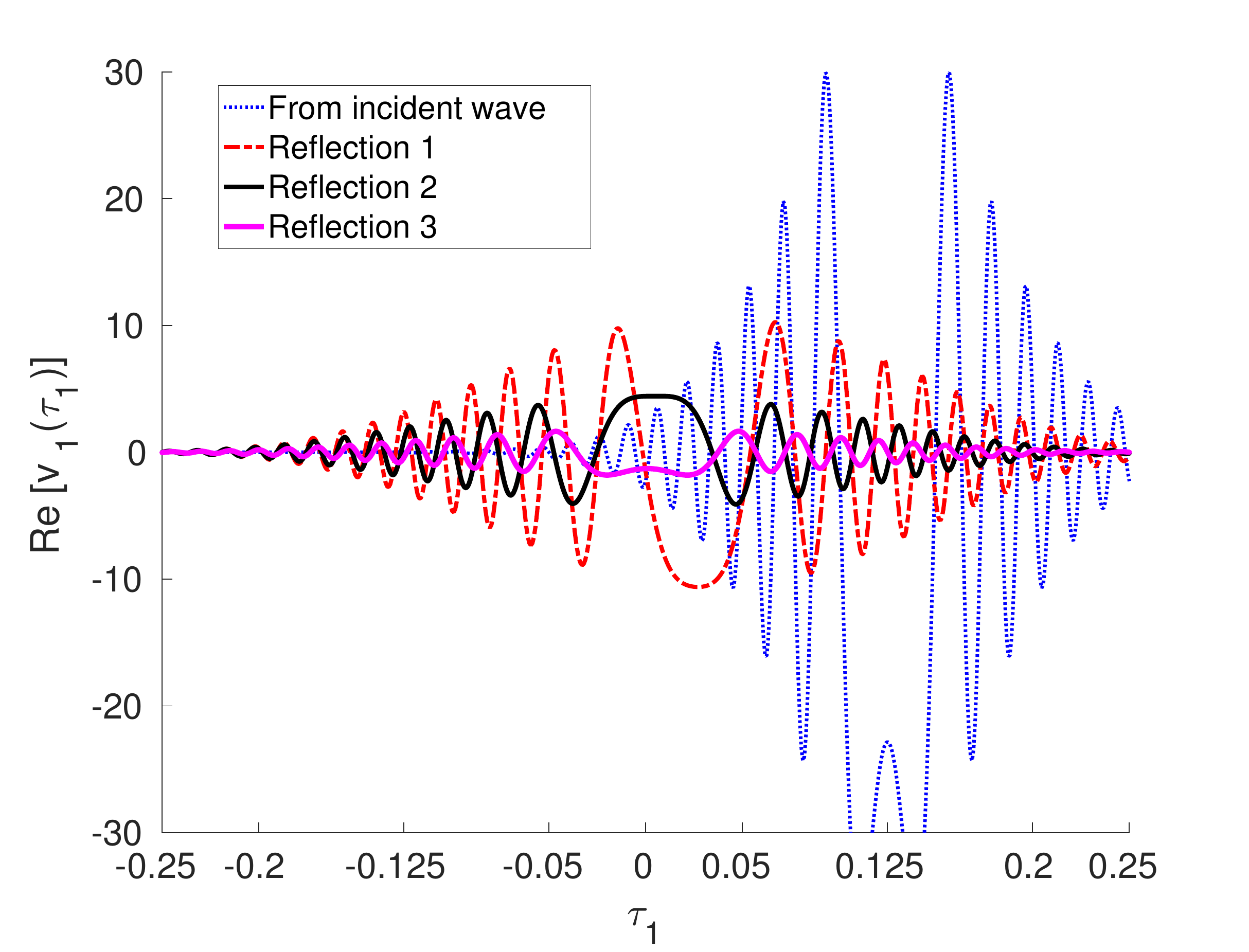} }
\caption[Multiple scattering iterations]{Consecutive densities in multiple scattering iterations at $k=2^7$ between the two disks shown in \cref{FtwoCircles}, with an incident plane wave from the left (left) and a point source at $(0,0)$ (right).}
\label{FmultScatIter}
\end{figure}

First, we illustrate in \cref{FmultScatIter} the convergence of the phase to a limiting function in consecutive iterations of the multiple scattering problem. We consider two circular obstacles again, which is similar to for example \cite[Fig. 1 \& 2]{geuzaine2005multiple} mentioned in \cref{Sintro}. The incident wave is a plane wave from the left (left panel) or a point source incidence (right panel). The density that is induced by the incident wave on the second disk is represented as a blue dotted line. In case of the plane wave, the density has a maximum and is least oscillatory near $\tau_2=-1/4$, which corresponds to the leftmost point on the disk. Already after the first reflection, the peak of the density is quite close to $\tau^*=0$. After two reflections, the peak is visually indistinguishable from the critical point $\tau^*$. Furthermore, the amplitude decreases by a certain factor after each reflection due to rays leaving the scene. In the right panel of the same figure, the results are repeated for a point source incident wave, and the conclusions are similar. These experiments were performed using submatrices from the numerical method outlined in the previous section.

Next, we turn to the computation of the limiting phase function in a range of $\tau$, rather than as a Taylor series. As outlined in \cref{Sgeometric}, we can compute $\chi$ from \cref{Echi}. We have applied six Newton-type iterations on \cref{Echi} with an initial guess given by $\chi(\tau) \approx (3-2\sqrt{2})\tau$, using the Chebfun package \cite{chebfun}. The resulting function satisfied \cref{Echi} with a residual close to machine precision. Furthermore, its Taylor series coefficients closely match our coefficients $a_i$. We have added a few iterations of $\chi$ as dashed purple lines in the left panel of \cref{FtwoCircles}. Recall that applying $\chi$ corresponds to tracing a ray back in time, and as such it contracts to the periodic orbit, as the figure illustrates.

\begin{figure}[h]
\centering
\subfloat{\includegraphics[width=0.45\hsize]{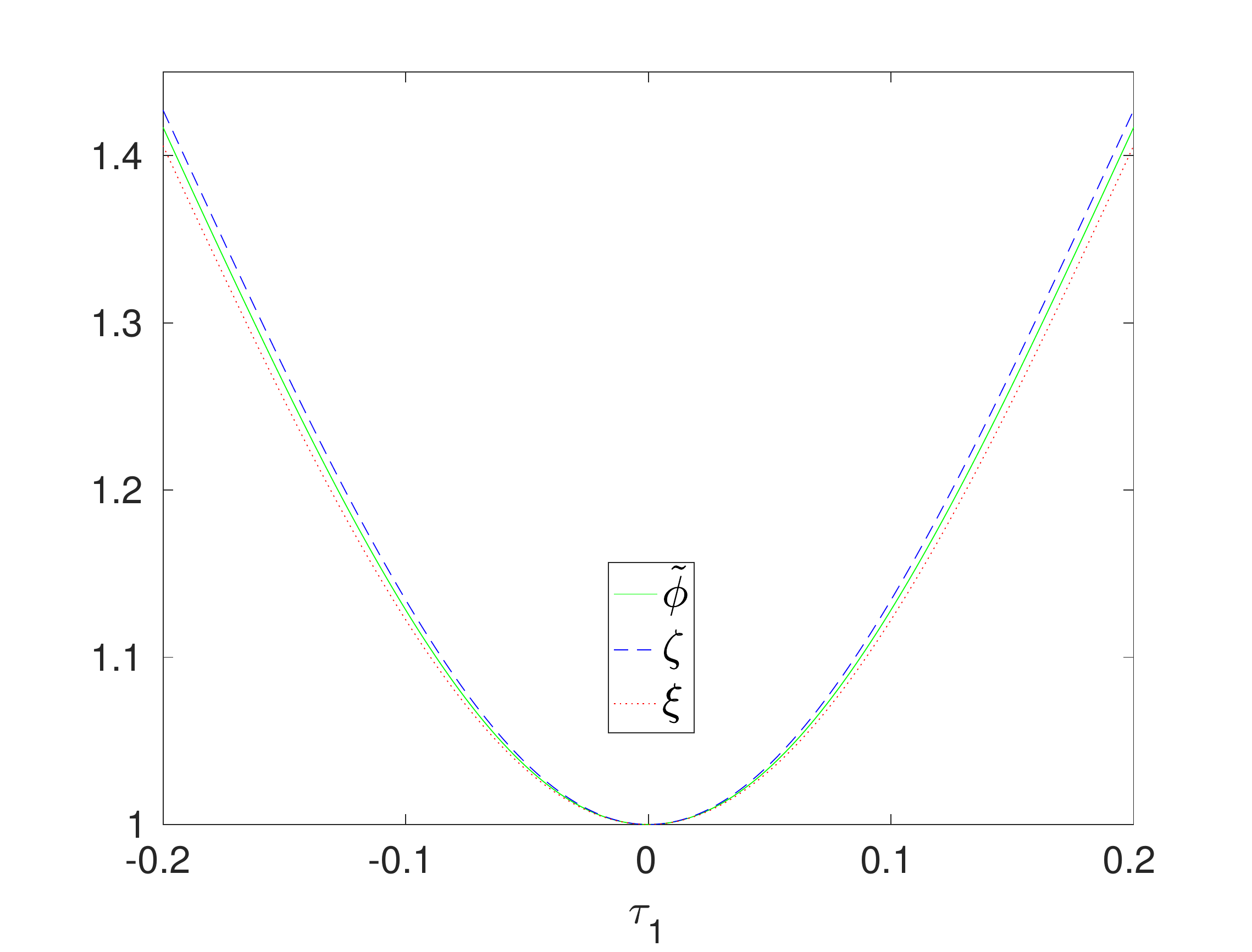} }
\subfloat{\includegraphics[width=0.54\hsize]{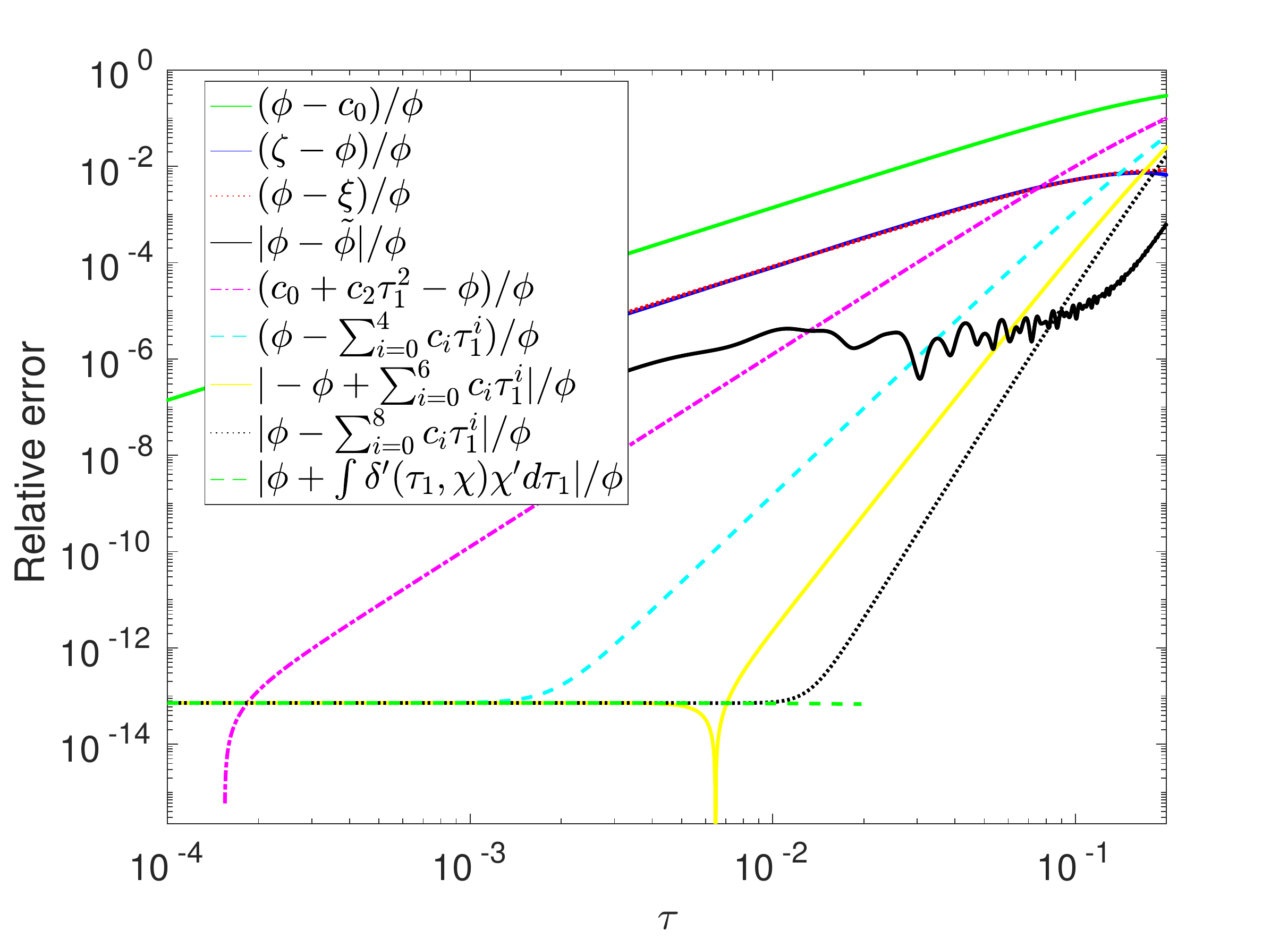} }
\caption[Phase and relative errors]{Left: comparison between the actual phase function $\phi(\tau)$, computed via \cref{EphiInf}, and the physical distances $\zeta$ and $\xi$ shown in \cref{FtwoCircles}. Right: relative errors of the various approximations to the phase considered in this paper. High accuracy is seen for the high-order Taylor series expansions for small $\tau$, i.e., close to the critical point.}
\label{Fphase}
\end{figure}

The right panel of \cref{Fphase} quantifies the convergence rate for a diverse set of approximations to the phase. Here, we have chosen a wavenumber $k=2^9$. We compute the phase using the geometric sum of distances given by \cref{EphiInf}. We end the summation after $r=10$ reflections, as the next term is close to machine precision in absolute value. The Taylor series with series expansion coefficients $c_i$ computed in \cref{Sexplicit} converges to this phase at a rate $\tau^T$ when taking $T$ terms. Note that the accuracy of $\tilde{\phi}$, recall \cref{Ephitildeprox}, is limited to $1e-6$ in this figure. This is due to the discretization error of the boundary element method, the composition of BEM submatrices, the computation of an eigenvector, approximating its argument and trying to undo the modulo $2\pi$ in it. The phase can indeed also be computed as \cref{EphiChi}, although in a smaller domain due to the inverse of $\chi$.

In the left panel of \cref{Fphase} we illustrate that the limiting phase seems well approximated by two physical distances $\zeta$ and $\xi$ that are easily computed and are illustrated in the left panel of \cref{FtwoCircles}. The value $\zeta$ is the distance from a point on $\Gamma_1$ to the point $\Gamma_2(0)$, i.e., the distance to the critical point on the other obstacle. The value $\xi$ is the distance from a point on $\Gamma_1$ to its closest point on $\Gamma_2$. Neither of these values is an exact expression of the phase, since they do not take into account reflections. However, we merely want to point out that both values are reasonable approximations, and they may prove to be useful as starting values for the solution of the non-linear system of equations \cref{EcircleSP}--\cref{EcirclePhase}.

\subsection{Two general obstacles}

\begin{figure}
\centering
\subfloat{\includegraphics[width=0.45\textwidth]{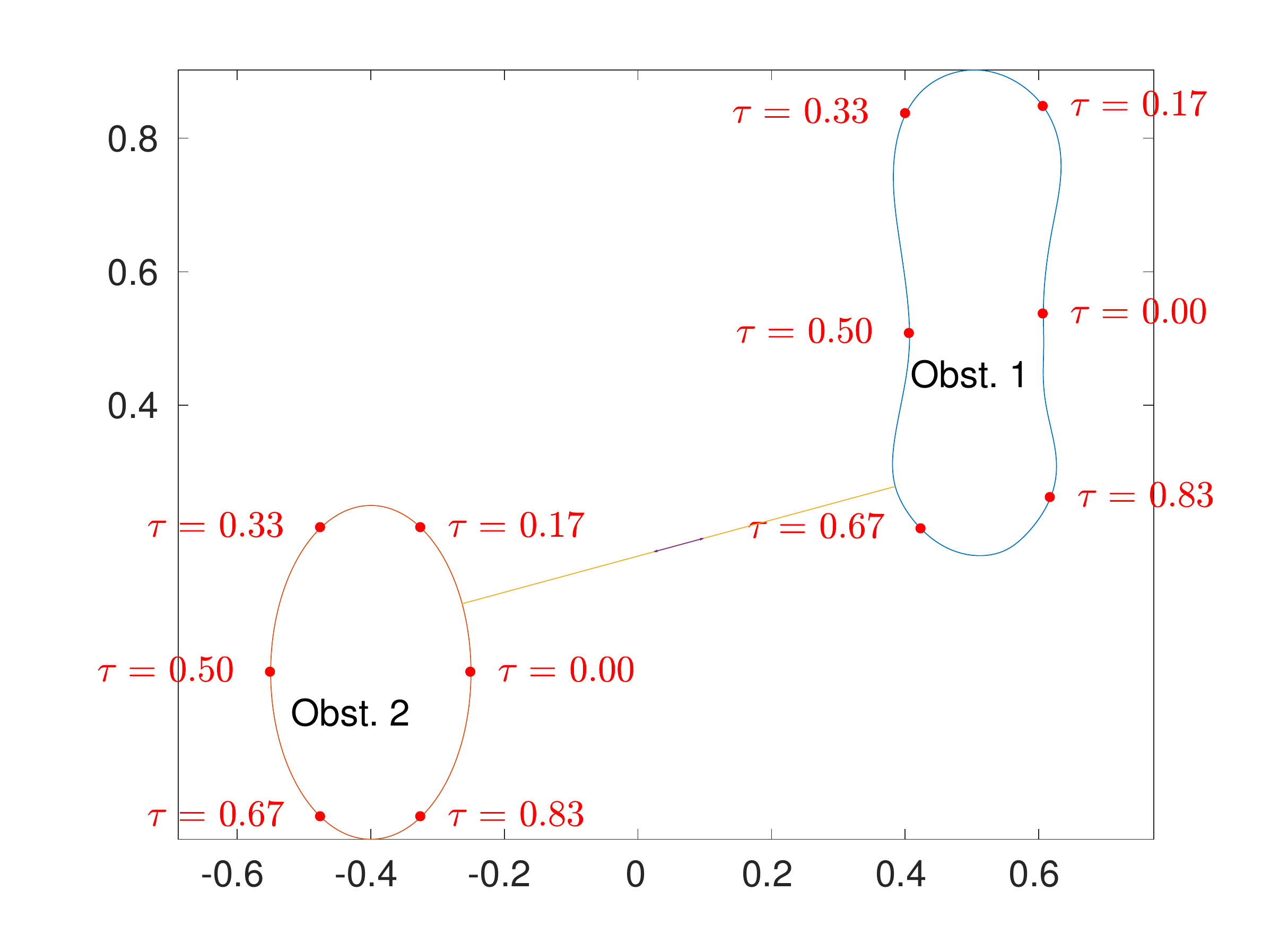} }
\subfloat{\includegraphics[width=0.45\textwidth]{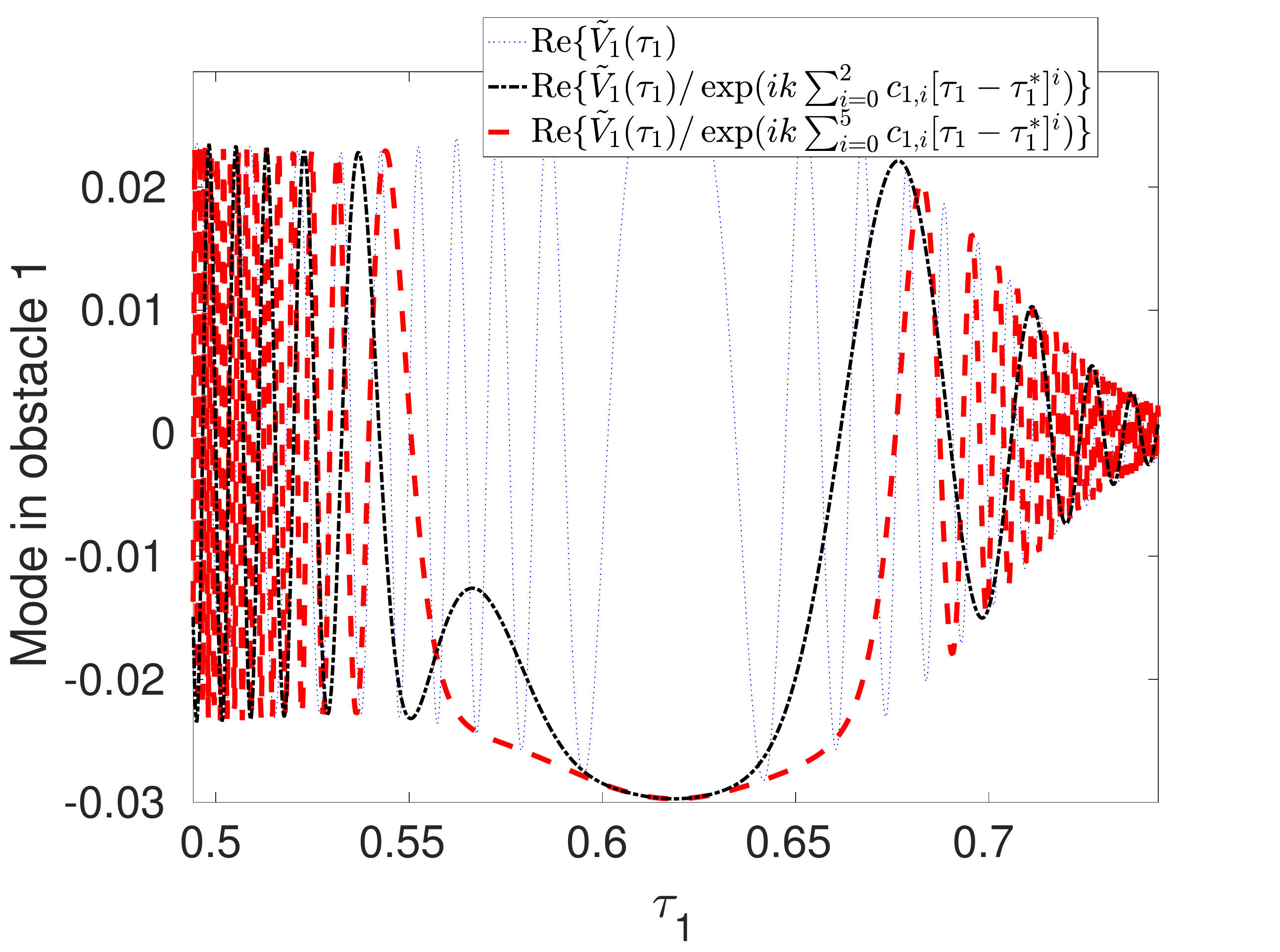} }
\caption[Periodic orbit for two obstacles]{Left: Periodic orbit between a near-convex obstacle and an ellipse, computed as a local minimum of \cref{eq:length}. Right: the real part of the corresponding mode $\tilde{V}_1(\tau_1)$ on the near-convex obstacle, computed for $k=512$.}
\label{FtwoObstMode}
\end{figure}

The scattering configuration shown in \cref{FtwoObstMode} does not exhibit the symmetry of the two disks. Moreover, the first obstacle is not convex. Still, we can consider the eigenmode corresponding to a periodic orbit connecting the two closest points on the two obstacles. In the right panel, we see that the mode is again locally nonoscillatory at $\tau_1^* \approx 0.62$ and its modulus still has a local maximum there. Our expressions for the Taylor series coefficients remain applicable. Factoring out a cubic or quintic approximation of the phase around $\tau_1^*$ does result in non-oscillatory behaviour of the eigenfunction near the critical point, as shown in the right panel of \cref{FtwoObstMode}. However, the quintic approximation does not seem to be much better than the cubic one.

\begin{figure}
\centering
\subfloat{\includegraphics[width=0.45\textwidth]{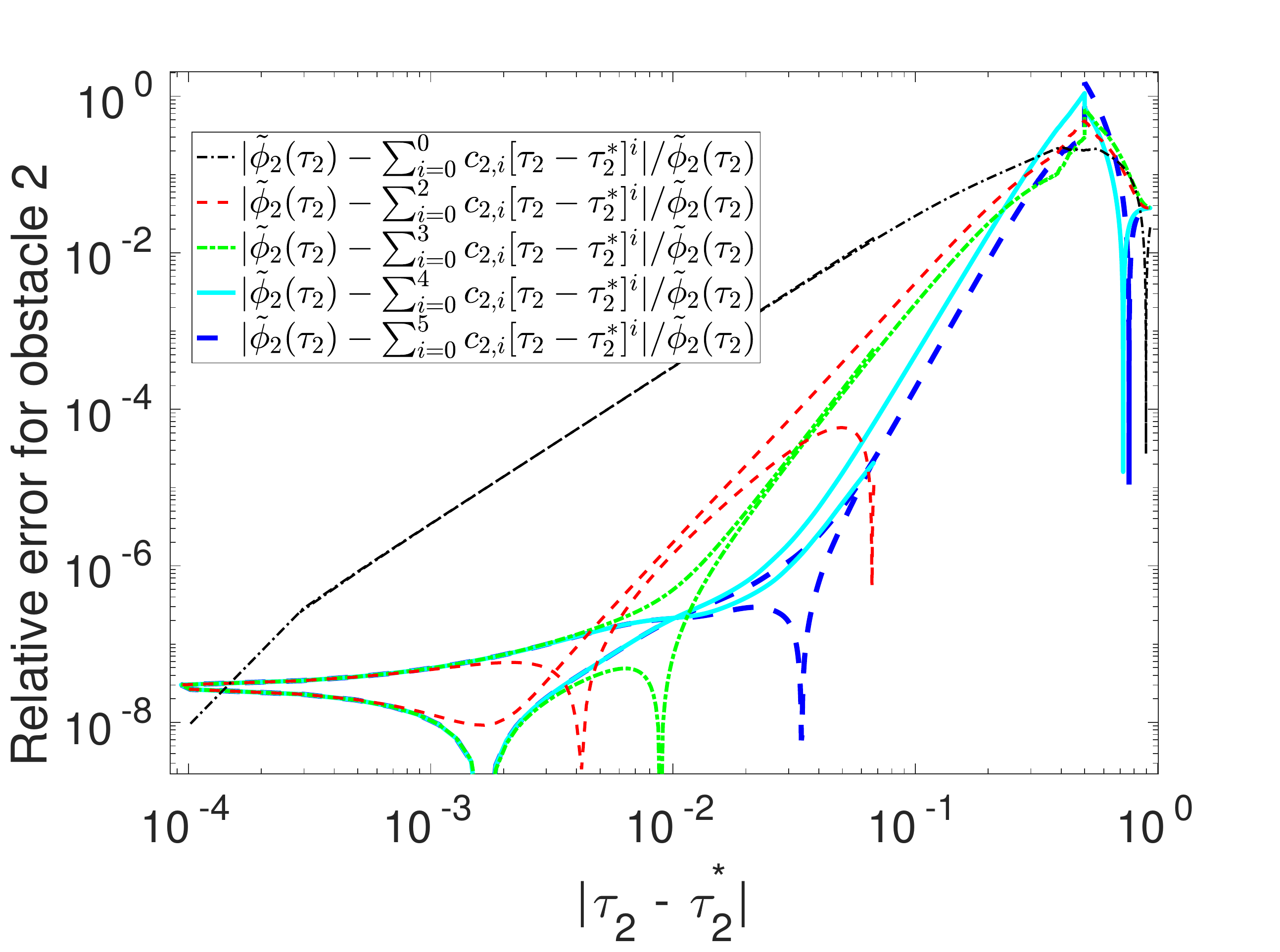}}
\subfloat{ \includegraphics[width=0.45\textwidth]{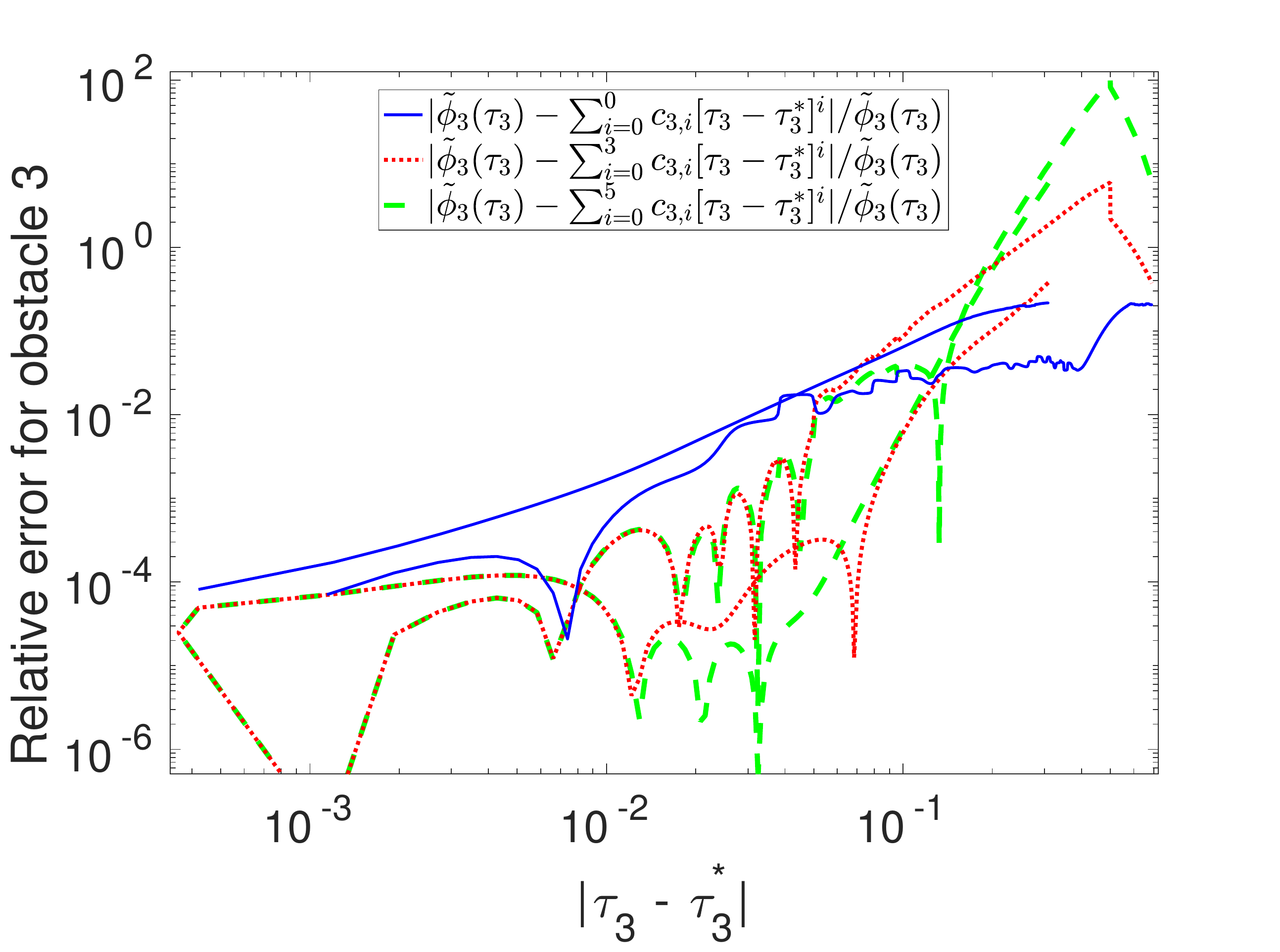} }
\caption[Convergence of the phase]{Comparison between the numerically computed phase $\tilde{\phi}$ calculated via \cref{Ephitildeprox}, and the Taylor series approximation with a varying number of terms, for scattering with two obstacles (left panel) and three obstacles (right panel).}
\label{FconvPhperOrbit}
\end{figure}

The accuracy is quantified in \cref{FconvPhperOrbit}. The numerical results for the two obstacles are shown in the left panel. Higher degree Taylor series approximations of the phase do lead to higher accuracy in a larger neighbourhood of the critical point. The accuracy in this figure is limited to $4e-8$ because of similar reasons as with the two disks. For relative errors well above this limit, we also see that the slope increases when adding more terms, as expected. So although less clear than the right part of \cref{Fphase}, \cref{FconvPhperOrbit} suggests we obtain the correct Taylor approximation of the phase.

We see a decline in \cref{FtwoObstMode} of the absolute value of the mode at the right of $\tau_1^* \approx 0.62$, which corresponds to the convex part at the bottom of $\Gamma_1$. More reflections are possible in the nonconvex part at the left of $\Gamma_1$, resulting in a slower decrease of $|V_1(\tau_1)|$ for $\tau_1$ decreasing from $0.62$. Consequently, one can now distinguish two different curves for each number of terms in \cref{FconvPhperOrbit}, where the one with the highest error corresponds to $\tau_1 < 0.62$.

\subsection{Three general obstacles} \label{SresThreeObsts}

\begin{figure}
\centering
\subfloat{\includegraphics[width=0.45\textwidth]{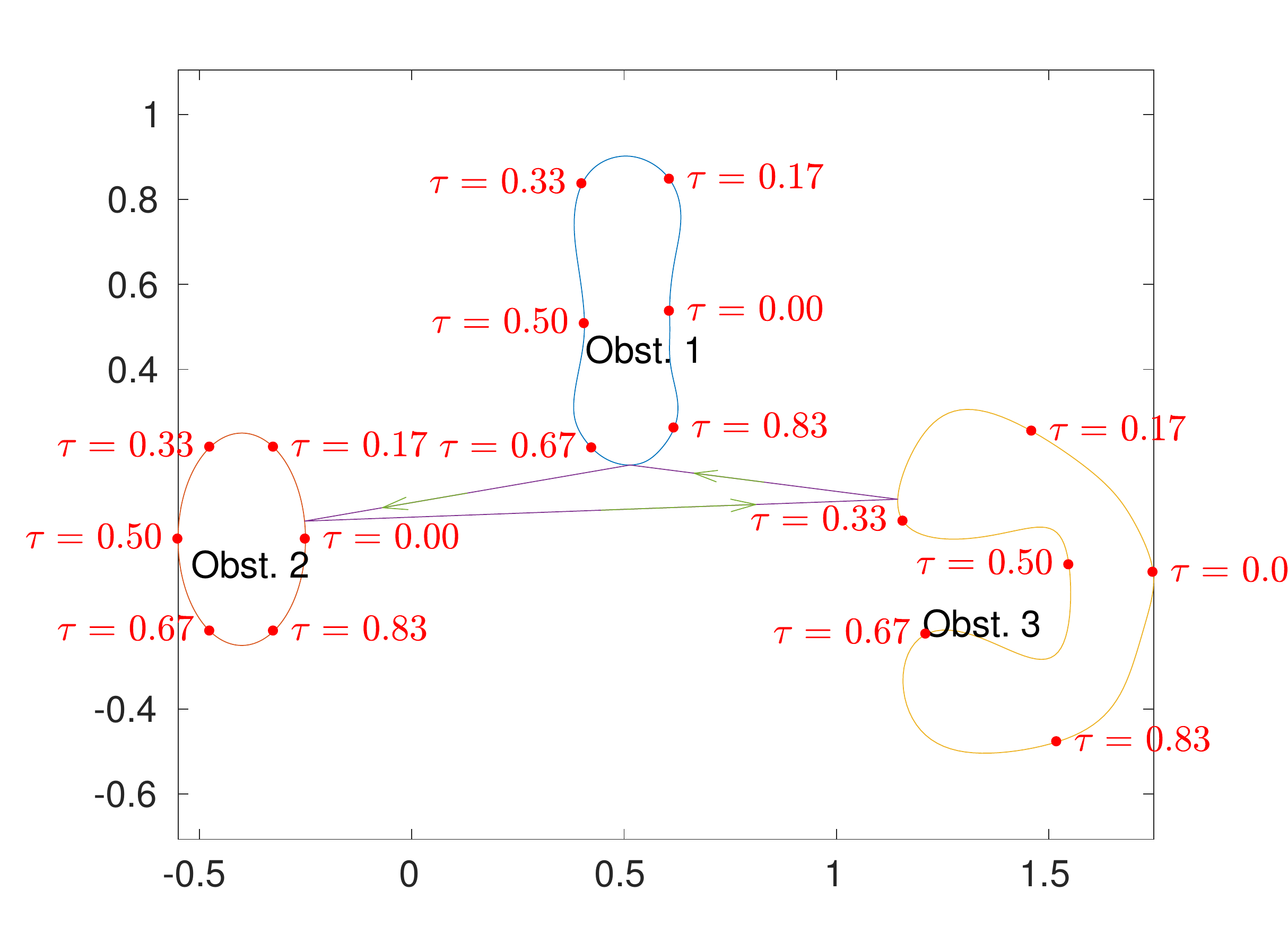}}
\subfloat{\includegraphics[width=0.45\textwidth]{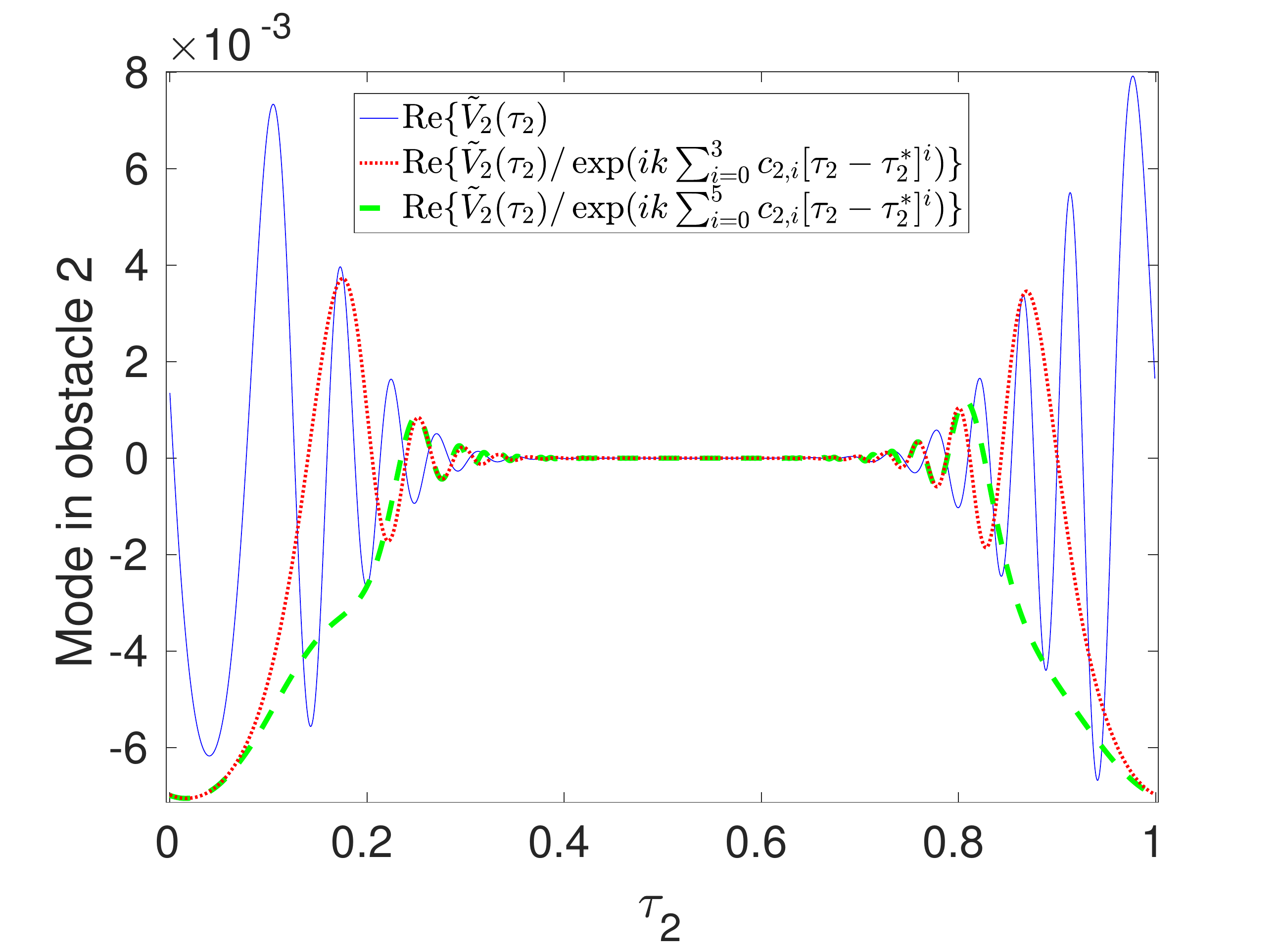}}
\caption[Periodic orbit for three obstacles]{Periodic orbit for near-convex obstacle, ellipse and near-inclusion (left) and real part of the corresponding mode $E_{j,1}$ on the ellipse for $k=128$ (right).}
\label{FthreeObstMode}
\end{figure}

Finally, we consider a configuration involving three obstacles ($J=3$), shown in \cref{FthreeObstMode}. Here, the critical point $\tau_3^*$ on the periodic orbit is in a convex part. However, the third obstacle is strongly non-convex, and has a near-inclusion with extremely complicated high-frequency behaviour in the range $0.33 < \tau_3 < 0.67$, $\tau_3^*$ is just outside this range and the periodic orbit is quite flat. Moreover, there is an additional $J=3$ periodic orbit near $\tau_3 = 0.67$, which may influence the one under investigation.

Still, even though the phase of the mode on the ellipse is influenced by the near-inclusion obstacle, the right part of \cref{FthreeObstMode} shows we obtain a good approximation of it. The error in the phase is quantified in the right panel of \cref{FconvPhperOrbit}, and it is clearly worse than in the previous example. Yet, somewhat surprisingly, even at the relatively small wavenumber $k=128$, high accuracy is reached near the critical point. The plots as a function of $|\tau_3 - \tau_3^*|$ now again have two clearly visible distinguished parts: the top one corresponds to $\tau_3 > \tau_3^*$ which is the inclusion part, while the bottom one corresponds to $\tau_{3} < \tau_3^*$, the convex part. The accuracy is clearly better in the convex part of the domain.

\section{Concluding remarks and future research} 

The asymptotic expressions in this paper are independent of the wavenumber and the incident wave. As a result, one can approximate the limiting phases of the densities on each obstacle in a periodic orbit at a cost that is independent of the wavenumber. That is a significant advantage. However, in any set of $J$ obstacles, there is potentially a very large number of periodic orbits. Indeed, there is a periodic orbit for any subset of the $J$ obstacles. This combinatorial cost of multiple scattering problems is inherent to the ray-tracing methodology: the cost of ray-tracing is independent of the wavenumber, but increases rapidly with an increasing number of obstacles.

Still, the results of this paper may be used to accelerate an iterative hybrid numerical-asymptotic method (such as the methods of Ecevit et al in \cite{2DEcevit,3DEcevit, EcevitOzen2017GBEMConvex,EcevitReitich2005Rate,boubendirEcevit2017acceleration}) as follows. A large cost of the ray-tracing is the long tail, i.e., a large number of iterations has to be performed before the difference between two consecutive iterations is small. Previous efforts have focused on accelerating the computation of the tail. However, the phases settle down quickly after only a few iterations. By extracting the asymptotic phase simultaneously for all densities on all obstacles at once, it may be possible to solve the tail at once as a coupled non-oscillatory problem. That is, we propose to perform a few iterations using phase-extraction in every step, followed by the solution of a coupled problem with the asymptotic phase extracted from each density. Since the consecutive iterations lead to the same phase in every cycle, the solution of the coupled problem is expected to be non-oscillatory. This is a topic of possible future research.

The cost of phase-extraction in the first few iterations can be further reduced if it is known how the phase of the incident wave evolves to the limiting phase we computed as a function of the number of reflections. A series expansion of $v_j^\text{smooth}(\tau)$ for the limiting density would even allow this without needing a mesh for $V_j(\tau)$ at all, from which also the magnitude of the associated eigenvalues can be derived.

\section*{Acknowledgments}
The authors would like to thank Samuel Groth and Niels S\o{}ndergaard for interesting and useful discussions on topics related to this paper. The authors were supported by FWO Flanders [projects G.0617.10, G.0641.11 and G.A004.14].

\bibliographystyle{abbrv}
\bibliography{perOrbit11arXiv.bbl}

\end{document}